\documentclass[11pt,a4paper]{article}

\usepackage{graphicx}
\usepackage{latexsym}
\usepackage{eepic}
\usepackage{amsfonts}
\makeatletter
\newdimen\mynormalparindent
\def\mymakefnmark{}%%%\hbox{${}^*$}}
\def\mymakefntext{\indent\mymakefnmark}
\long\def\myfootnotetext#1{\insert\footins{%
  \normalfont\footnotesize
  \interlinepenalty\interfootnotelinepenalty
  \splittopskip\footnotesep \splitmaxdepth \dp\strutbox
  \floatingpenalty\@MM \hsize\columnwidth
  \@parboxrestore \parindent\mynormalparindent \sloppy
  \mymakefntext{\rule\z@\footnotesep\ignorespaces#1\unskip\strut\par}}}

\makeatother

\newtheorem{basic}{Basic}[section]

\newtheorem{lem}[basic]{Lemma}
\newtheorem{propos}[basic]{Proposition}
\newtheorem{thm}[basic]{Theorem}
\newtheorem{exa}[basic]{Example}

\newtheorem{cor}[basic]{Corollary}

\newtheorem{conjec}[basic]{Conjecture}

\newcommand{\bdm}{\begin{displaymath}}
\newcommand{\edm}{\end{displaymath}}
\newcommand{\be}{\begin{equation}}
\newcommand{\ee}{\end{equation}}
\newcommand{\ep}{\vspace{-3mm}\hfill\mbox{$\Box$}\\}
\newcommand{\ul}{\underline}

\newcommand{\R}{\mathbb{R}}

\def\theequation{\thesection.\arabic{equation}}

\usepackage{color}
\definecolor{purple}{rgb}{0.5,0,0.3}

\begin{document}

\begin{center}
{\bf\Large Subtracting a best rank-1 approximation may
increase tensor rank}
\footnote{A 5-page summary of this paper has been accepted for the EUSIPCO
2009 conference and will appear in its proceedings.}\\ \vspace{2cm}

Alwin Stegeman\footnote{Corresponding author. A. Stegeman is with the
Heijmans Institute for Psychological Research, University of Groningen,
Grote Kruisstraat 2/1, 9712 TS Groningen, The Netherlands, phone: +31 50
363 6193, fax: +31 50 363 6304, email: a.w.stegeman@rug.nl, URL:
http://www.gmw.rug.nl/$\sim$stegeman. Research supported by the Dutch
Organisation for Scientific Research (NWO), VENI grant 451-04-102 and VIDI
grant 452-08-001.}
and Pierre Comon\footnote{P. Comon is with the Laboratoire I3S, UMR6070
CNRS, UNS, 2000, route des Lucioles, BP 121, 06903 Sophia Antipolis Cedex,
France, phone: + 33 4 92 94 27 17, fax: +33 4 92 94 28 98, email:
pcomon@unice.fr. Research supported by contract ANR-06-BLAN-0074 ``Decotes".
}\\ \vspace{1cm} \today \vspace{1cm} \end{center}

\begin{abstract}
\noindent It has been shown that a best rank-$R$ approximation of an
order-$k$ tensor may not exist when $R\ge 2$ and $k\ge 3$. This poses a
serious problem to data analysts using tensor decompositions. It has been
observed numerically that, generally, this issue cannot be
solved by consecutively computing and subtracting best rank-1
approximations. The reason for this is that subtracting a best rank-1
approximation generally does not decrease tensor rank. In this paper, we
provide a mathematical treatment of this property for real-valued $2\times
2\times 2$ tensors, with symmetric tensors as a special case. Regardless of
the symmetry, we show that for generic $2\times 2\times 2$ tensors (which
have rank 2 or 3), subtracting a best rank-1 approximation results in a
tensor that has rank 3 and lies on the boundary between the rank-2 and
rank-3 sets. Hence, for a typical tensor of rank 2, subtracting a best
rank-1 approximation {\em increases} the tensor rank. \end{abstract}

\newpage
\section{Introduction}
\setcounter{equation}{0}
Tensors of order $d$ are defined on the outer product of $d$
linear spaces, ${\cal S}_\ell$, $1\le\ell\le d$. Once bases of
spaces ${\cal S}_\ell$ are fixed, they can be represented by
$d$-way arrays. For simplicity, tensors are usually assimilated with their
array representation. We assume throughout the following notation:
underscored bold uppercase for tensors e.g. $\ul{\bf X}$, bold uppercase
for matrices e.g. ${\bf T}$, bold lowercase for vectors e.g. ${\bf a}$,
calligraphic for sets e.g. ${\cal S}$, and plain font for scalars e.g.
$X_{ijk}$, $T_{ij}$ or $a_i$. In this paper, we consider only $3$rd
order tensors. The three spaces of a $3$rd order tensor are also referred
to as the three ``modes".

Let $\ul{\bf X}$ be a $3$rd order tensor defined on the tensor product
${\cal S}_1\otimes {\cal S}_2 \otimes {\cal S}_3$. If a
change of bases is performed in the spaces
${\cal S}_1,{\cal S}_2,{\cal S}_3$ by invertible matrices
${\bf S},{\bf T},{\bf U}$, then the tensor representation $\ul{\bf X}$ is
transformed into \begin{equation} \label{multilin-eq} \widetilde{\ul{\bf
X}} = ({\bf S}, {\bf T}, {\bf U})\cdot \ul{\bf X}\,, \end{equation}

\noindent whose coordinates are given by $\widetilde{X}_{ijk} = \sum_{pqr}
S_{ip} \, T_{jq} \, U_{kr} \, X_{pqr}$. This is known as the
\emph{multilinearity property} enjoyed by tensors. Matrices, which can be
associated with linear operators, are tensors of order 2.
The multilinear transformation (\ref{multilin-eq}) is also denoted as
\be
\label{eq-multicontr}
\widetilde{\ul{\bf X}} = \ul{\bf X} \bullet_1 {\bf S} \bullet_2 {\bf
T} \bullet_3 {\bf U}\,, \ee

\noindent where $\bullet_{\ell}$ denotes the multiplication (or
contraction) operator in the $\ell$th mode. Note that the matrix
multiplication ${\bf STU}^T$ can be denoted as ${\bf T}\bullet_1 {\bf
S}\bullet_2 {\bf U}={\bf T}\bullet_2 {\bf U}\bullet_1 {\bf S}$. For
two contractions with matrices in the same mode, we have the rule $\ul{\bf
X}\bullet_{\ell} {\bf T}\bullet_{\ell} {\bf S}=\ul{\bf X}\bullet_{\ell}
({\bf ST})$, see e.g. \cite[section 2]{LMVrank1}.

The rank of a tensor $\ul{\bf X}$ is defined as the smallest number
of outer product tensors whose sum equals $\ul{\bf X}$, i.e. the smallest
$R$ such that
\be
\ul{\bf X}=\sum_{r=1}^R {\bf a}_r \otimes {\bf b}_r \otimes {\bf c}_r\,.
\ee

\noindent Hence a rank-1 tensor $\ul{\bf X}$ is the outer product of
vectors ${\bf a},{\bf b},{\bf c}$ and has entries $X_{ijk}=a_ib_jc_k$.
The decomposition of a tensor into a sum of outer products of vectors and
the corresponding notion of tensor rank were first introduced and studied
by \cite{H1} \cite{H2}.

Tensors play a wider and wider role in numerous application areas including
blind source separation techniques for Telecommunications
\cite{SGB} \cite{SBG2} \cite{Comon} \cite{LC} \cite{AlmFM}, Arithmetic
Complexity \cite{K} \cite{Stras} \cite{BCS} \cite{SS}, or Data
Analysis \cite{Har} \cite{CC} \cite{SBG} \cite{Kro}. In some applications,
tensors may be symmetric only in some modes, or may not
be symmetric nor have equal dimensions. In most applications, the
decomposition of a tensor into a sum of rank-1 terms is relevant, since
tensors entering the models to fit have a reduced rank. For
example, such a tensor decomposition describes the basic structure of
fourth-order cumulants of multivariate data on which a lot of algebraic
methods for Independent Component Analysis are based \cite{pierre}
\cite{LMV2} \cite{Hyv}. For an overview of applications of tensor
decompositions, refer to \cite{KolBad}.

An important advantage of using tensor decompositions of order 3 and
higher, is that the decomposition is rotationally unique under mild
conditions \cite{K} \cite{SS}. This is not the case for most matrix
decompositions, e.g. Principal Component Analysis. However, the
manipulation of tensors remains difficult, because of major differences
between their properties when we go from second order to higher orders. We
mention the following: $(i)$ tensor rank often exceeds dimensions,
$(ii)$ tensor rank can be different over the real and complex fields,
$(iii)$ maximal tensor rank is not generic, and is still unknown in
general, $(iv)$ generic tensor rank may not have a single value over the
real field, $(v)$ computing the rank of a tensor is very difficult, $(vi)$
a tensor may not have a best rank-$R$ approximation for $R\ge 2$. For
$(i)$-$(v)$, see e.g. \cite{Krus89} \cite{DSL} \cite{CGLM}. For $(iv)$,
see e.g. \cite{TBK} \cite{TBSR}. For $(vi)$, see e.g. \cite{Ste}
\cite{Ste3} \cite{Ste2} \cite{DSL} \cite{KDS} \cite{SDL}. A discussion
specifically focussed on symmetric tensors can be found in \cite{CGLM}.

In \cite{DSL} it is shown that $(vi)$ holds on a set of positive measure.
It is recalled in \cite{CGLM} and \cite{DSL} that any tensor has a best
rank-1 approximation. However, it has been observed numerically in
\cite[section 7]{KofReg} that a best or "good" rank-$R$ approximation
cannot be obtained by consecutively computing and substracting $R$ best
rank-1 approximations. The reason for this is that subtracting a best
rank-1 approximation generally does not decrease tensor rank. Hence, the
deflation technique practiced for matrices (via the Singular Value
Decomposition) cannot generally be extended to higher-order tensors. A
special case where this deflation technique works is when the tensor is
diagonalizable by orthonormal multilinear transformation; see \cite[section
7]{KofReg}.

In this paper, we provide a mathematical treatment of the (in)validity of
a rank-1 deflation procedure for higher-order tensors. We consider $2\times
2\times 2$ tensors over the real field. For such a tensor $\ul{\bf X}$, let
the frontal slabs be denoted as ${\bf X}_1$ and ${\bf X}_2$. Our main
result is for generic tensors $\ul{\bf X}$, which have rank 2 if ${\bf
X}_2{\bf X}_1^{-1}$ has distinct real eigenvalues, and rank 3 if ${\bf
X}_2{\bf X}_1^{-1}$ has complex eigenvalues. We show that for generic
$\ul{\bf X}$, subtraction of a best rank-1 approximation $\ul{\bf Y}$
yields a tensor $\ul{\bf Z}=\ul{\bf X}-\ul{\bf Y}$ of rank 3. Hence, for
a typical $\ul{\bf X}$ of rank 3 this does not affect the rank, and for
a typical $\ul{\bf X}$ of rank 2 this has {\em increased} the rank. In
fact, we show that $\ul{\bf Z}$ lies on the boundary between the rank-2
and rank-3 sets, i.e. ${\bf Z}_2{\bf Z}_1^{-1}$ has identical real
eigenvalues. The result that subtraction of a best rank-1 approximation
yields identical eigenvalues is new and expands the knowledge of the
topology of tensor rank. Also, we show that the same result holds for
symmetric $2\times 2\times 2$ tensors. Based on numerical experiments we
conjecture that the results can be extended to $p\times p\times 2$ tensors
over the real field.

The above contributions are new to the
literature on best rank-1 approximation of higher-order tensors. The latter
includes best rank-1 approximation algorithms \cite{LMVrank1} \cite{Zhang}
\cite{KofReg}, conditions under which the best rank-1 approximation is
equal to the best symmetric rank-1 approximation \cite{Ni}, and a relation
between the best symmetric rank-1 approximation and the notions of
eigenvalues and eigenvectors of a symmetric tensor \cite{Com98} \cite{Qi}.

This paper is organized as follows. In Section 2, we introduce the best
rank-1 approximation problem for $3$rd order tensors, and state first order
conditions for the optimal solution. Next, we consider $2\times 2\times 2$
tensors. Section 3 contains rank criteria and orbits for $2\times
2\times 2$ tensors. In Section 4, we present examples and general results
for subtraction of a best rank-1 approximation from a $2\times 2\times 2$
tensor. In Section 5, 6 and 7, we discuss the special case of symmetric
tensors. Section 5 provides first order conditions for the best symmetric
rank-1 approximation of a symmetric $3$rd order tensor. Section 6 contains
rank criteria and orbits of symmetric $2\times 2\times 2$ tensors. These
results are used in Section 7, when studying the subtraction of a best
symmetric rank-1 approximation from a symmetric $2\times 2\times 2$ tensor.
Section 8 contains a discussion of our results. The proofs of our
main results are contained in appendices.

\section{Best rank-1 approximation}
\setcounter{equation}{0}
We consider the problem of finding a best rank-1 approximation to a given
$3$rd order tensor $\ul{\bf X}\in\R^{d_1\times d_2\times d_3}$, i.e.
\be \label{eq-prob} {\rm min}_{{\bf x}\in\R^{d_1},{\bf y}\in\R^{d_2},{\bf
z}\in\R^{d_3}}\;||\ul{\bf X}-{\bf x}\otimes{\bf y}\otimes{\bf z}||^2\,, \ee

\noindent where $||\cdot||$ denotes the Frobenius norm, i.e. $||\ul{\bf
X}||^2 = \sum_{ijk} |X_{ijk}|^2$. Since the set of rank-1 tensors is
closed, problem (\ref{eq-prob}) is guaranteed to have an optimal solution
\cite[proposition 4.2]{DSL}. Note that the vectors ${\bf x},{\bf y},{\bf
z}$ of the rank-1 tensor $({\bf x}\otimes{\bf y}\otimes{\bf z})$ are
determined up to scaling. One could impose two of the vectors to be unit
norm.

The criterion of (\ref{eq-prob}) can be written as
\be
\label{eq-crit}
\Psi=||\ul{\bf X}||^2 -2\,\ul{\bf X}\bullet_1 {\bf x}^T\bullet_2 {\bf y}^T
\bullet_3 {\bf z}^T + ||{\bf x}||^2 ||{\bf y}||^2 ||{\bf z}||^2\,.
\ee

\noindent When the gradients of $\Psi$ with respect to the vectors ${\bf
x},{\bf y},{\bf z}$ are set to zero, we obtain the following equations:
\be
\label{eq-grad}
{\bf x}  = \frac{\ul{\bf X} \bullet_2 {\bf
y}^T \bullet_3 {\bf z}^T}{||{\bf y}||^2 ||{\bf z}||^2}\,, \quad\quad
{\bf y}  = \frac{\ul{\bf X} \bullet_1 {\bf x}^T
\bullet_3 {\bf z}^T}{||{\bf x}||^2 ||{\bf z}||^2}\,, \quad\quad
{\bf z} = \frac{\ul{\bf X} \bullet_1 {\bf x}^T \bullet_2
{\bf y}^T}{||{\bf x}||^2 ||{\bf y}||^2}\,. \ee

\noindent Substituting
\be
\label{eq-xopt}
{\bf x}= \frac{\ul{\bf X}\bullet_2 {\bf y}^T\bullet_3 {\bf z}^T}
{||{\bf y}||^{2} ||{\bf z}||^{2}}\,,
\ee

\noindent into the last two equations of (\ref{eq-grad}), we obtain
\be
\label{eq-gradsub}
(\ul{\bf X}\bullet_3 {\bf z}^T) \bullet_1 (\ul{\bf X}\bullet_3 {\bf z}^T)
\bullet_2 {\bf y}^T = \lambda\, {\bf y}\,,\quad\quad\quad\quad
(\ul{\bf X}\bullet_2 {\bf y}^T) \bullet_1 (\ul{\bf X}\bullet_2 {\bf y}^T)
\bullet_3 {\bf z}^T = \mu\, {\bf z}\,,
\ee

\noindent where $\lambda=||{\bf x}||^2 ||{\bf y}||^2 ||{\bf z}||^4$ and
$\mu=||{\bf x}||^2 ||{\bf y}||^4 ||{\bf z}||^2$. Hence, ${\bf y}$ is an
eigenvector of the matrix $(\ul{\bf X}\bullet_3 {\bf z}^T) \bullet_1
(\ul{\bf X}\bullet_3 {\bf z}^T)$ and ${\bf z}$ is an eigenvector of the
matrix $(\ul{\bf X}\bullet_2 {\bf y}^T) \bullet_1 (\ul{\bf X}\bullet_2 {\bf
y}^T)$.

Substituting (\ref{eq-xopt}) into the criterion
(\ref{eq-crit}) yields \be \label{eq-critsub} \Psi = ||\ul{\bf X}||^2
-\frac{(\ul{\bf X}\bullet_2 {\bf y}^T\bullet_3 {\bf z}^T)\bullet_1 (\ul{\bf
X}\bullet_2 {\bf y}^T\bullet_3 {\bf z}^T)} {||{\bf y}||^2\,||{\bf z}||^2}=
||\ul{\bf X}||^2 -\frac{||\ul{\bf X}\bullet_2 {\bf y}^T\bullet_3 {\bf
z}^T||^2} {||{\bf y}||^2\,||{\bf z}||^2}\,. \ee

\noindent Hence, a best rank-1 approximation $({\bf x}\otimes {\bf y}\otimes
{\bf z})$ of $\ul{\bf X}$ is found by minimizing (\ref{eq-critsub}) over
$({\bf y},{\bf z})$ and obtaining ${\bf x}$ as (\ref{eq-xopt}). The
stationary points $({\bf y},{\bf z})$ are given by (\ref{eq-gradsub}),
which can also be written as
\begin{eqnarray}
\label{eq-gradsub2y}
(\ul{\bf X}\bullet_2 {\bf y}^T\bullet_3 {\bf z}^T)\bullet_1
(\ul{\bf X}\bullet_3 {\bf z}^T) &=&  \frac{||\ul{\bf X}\bullet_2 {\bf
y}^T\bullet_3 {\bf z}^T||^2}{||{\bf y}||^2}\;{\bf y}\,, \\[2mm]
\label{eq-gradsub2z}
(\ul{\bf X}\bullet_2 {\bf y}^T\bullet_3 {\bf z}^T)\bullet_1 (\ul{\bf
X}\bullet_2 {\bf y}^T) &=&  \frac{||\ul{\bf X}\bullet_2 {\bf y}^T\bullet_3
{\bf z}^T||^2}{||{\bf z}||^2}\;{\bf z}\,.
\end{eqnarray}

\noindent Next, we consider transformations of the best rank-1
approximation. The following well-known result states that a best
rank-1 approximation is preserved under orthonormal multilinear
transformation.

\begin{lem}
\label{p-orth}
Let ${\bf S},{\bf T},{\bf U}$ be orthonormal matrices.
If a tensor $\ul{\bf X}$ admits $\ul{\bf Y}$ as a best rank-$1$
approximation, then $({\bf S},{\bf T},{\bf U})\cdot\ul{\bf Y}$ is a best
rank-$1$ approximation of $({\bf S},{\bf T},{\bf U})\cdot\ul{\bf X}$.
\end{lem}

\noindent {\bf Proof.} Let $\ul{\bf Y}={\bf x}\otimes{\bf y}\otimes{\bf
z}$ be a best rank-1 approximation of $\ul{\bf X}$, and let
$\ul{\widetilde{\bf X}}=({\bf S},{\bf T},{\bf U})\cdot\ul{\bf X}$. Since
orthonormal transforms leave the Frobenius norm invariant, we obtain the
following analogue of (\ref{eq-crit}):
\begin{eqnarray}
\label{eq-critorth}
~\hspace{-1cm}||\ul{\widetilde{\bf X}}-\tilde{\bf x}\otimes\tilde{\bf
y}\otimes\tilde{\bf z}||^2 &=& ||\ul{\widetilde{\bf
X}}||^2-2\,\ul{\widetilde{\bf X}}\bullet_1 \tilde{\bf
x}^T\bullet_2\tilde{\bf y}^T\bullet_3\tilde{\bf z}^T+ ||\tilde{\bf
x}||^2||\tilde{\bf y}||^2||\tilde{\bf z}||^2 \nonumber \\[2mm] &=&
||\ul{\bf X}||^2 -2\,\ul{\bf X}\bullet_1 (\tilde{\bf x}^T{\bf S})\bullet_2
(\tilde{\bf y}^T{\bf T}) \bullet_3 (\tilde{\bf z}^T{\bf U}) + ||{\bf
S}^T\tilde{\bf x}||^2 ||{\bf T}^T\tilde{\bf y}||^2 ||{\bf U}^T\tilde{\bf
z}||^2\,. \end{eqnarray}

\noindent Hence, since $({\bf x},{\bf y},{\bf z})$ is a minimizer of
(\ref{eq-crit}), it follows that $({\bf Sx},{\bf Ty},{\bf Uz})$ is a
minimizer of (\ref{eq-critorth}). In other words, a best rank-1
approximation of $\ul{\widetilde{\bf X}}$ is given by $({\bf Sx}\otimes{\bf
Ty}\otimes{\bf Uz})$. This completes the proof.\ep

\noindent As we will see later, most tensors have multiple locally best
rank-1 approximations, with one of them being better than the others
(i.e., a unique global best rank-1 approximation). Our final result in this
section states a condition under which there exist infinitely many best
(global) rank-1 approximations.

\begin{propos}
\label{p-infbest}
Let $\ul{\bf X}$ be such that the matrix $(\ul{\bf X}\bullet_3{\bf z}^T)$ is
orthogonal for any vector ${\bf z}$, and $(\ul{\bf X}\bullet_2{\bf y}^T)$
is orthogonal for any vector ${\bf y}$. Then $\ul{\bf X}$ has infinitely
many best rank-$1$ approximations.
\end{propos}

\noindent {\bf Proof.} The proof follows from equation
(\ref{eq-gradsub}) for the stationary points $({\bf y},{\bf z})$. The
conditions of the proposition imply that the matrices
$(\ul{\bf X}\bullet_3 {\bf z}^T) \bullet_1 (\ul{\bf X}\bullet_3 {\bf z}^T)$
and $(\ul{\bf X}\bullet_2 {\bf y}^T) \bullet_1 (\ul{\bf X}\bullet_2 {\bf
y}^T)$ are proportional to the identity matrix for any ${\bf y}$ and ${\bf
z}$. Therefore, any vector is an eigenvector of these matrices, and
(\ref{eq-gradsub}) holds for any $({\bf y},{\bf z})$.

Since any $({\bf y},{\bf z})$ is a stationary point of criterion
(\ref{eq-critsub}), it follows that the latter is constant. We conclude
that any $({\bf x}\otimes{\bf y}\otimes{\bf z})$ with ${\bf x}$ as in
(\ref{eq-xopt}), is a best rank-1 approximation of $\ul{\bf X}$.\ep

\noindent Below is a $2\times 2\times 2$ example satisfying the conditions
of Proposition~\ref{p-infbest}. We denote a tensor $\ul{\bf X}$ with two
slabs ${\bf X}_1$ and ${\bf X}_2$ as $[{\bf X}_1\,|\,{\bf X}_2]$.

\begin{exa}
\label{ex-KHL}
{\rm Let
\be
\label{eq-KHL}
\ul{\bf X} = \left[\begin{array}{cc|cc}1 & 0 & 0 & -1 \\ 0 & 1 & 1
& 0\end{array} \right]\,.
\ee

\noindent Then for any choice of nonzero vector ${\bf z}$,
the matrix $(\ul{\bf X} \bullet_3 {\bf z}^T)$, obtained by linear
combination of the above two matrix slices, is orthogonal. Also, for any
nonzero vector ${\bf y}$, the matrix $(\ul{\bf X} \bullet_2 {\bf y}^T)$ is
orthogonal. Hence, $\ul{\bf X}$ has infinitely many rank-1 approximations.
One can verify that each rank-1 approximation $({\bf x}\otimes{\bf
y}\otimes{\bf z})$ with ${\bf x}$ as in (\ref{eq-xopt}), satsifies
$||\ul{\bf X}-{\bf x}\otimes{\bf y}\otimes{\bf z}||^2=3$.}\ep \end{exa}

\noindent The tensor (\ref{eq-KHL}) has rank 3 and is studied in
\cite{BKL} where it is shown that it has no best rank-2
approximation, the infimum of $||\ul{\bf X}-\ul{\bf Y}||^2$ over $\ul{\bf
Y}$ of rank at most 2 being 1. A more general result is obtained in
\cite{DSL} where it is shown that no $2\times 2\times 2$ tensor
of rank 3 has a best rank-2 approximation. In \cite{Ste} it is shown that
any sequence of rank-2 approximations $\ul{\bf Y}^{(n)}$ for which
$||\ul{\bf X}-\ul{\bf Y}^{(n)}||^2$ converges to the infimum of 1, features
diverging components.

\section{Rank criteria and orbits of $2\times 2\times 2$ tensors}
\setcounter{equation}{0}
It was shown in \cite[section 7]{DSL} that $2\times
2\times 2$ tensors (over the real field) can be transformed by invertible
multilinear multiplications (\ref{multilin-eq}) into eight distinct
canonical forms. This partitions the space $\R^{2\times 2\times 2}$ into
eight distinct orbits under the action of invertible transformations of a
tensor ``from the three sides".

Before the eight orbits are introduced, we define some concepts. A mode-$n$
vector of a $d_1\times d_2\times d_3$ tensor is an $d_n\times 1$ vector
obtained from the tensor by varying the $n$-th index and keeping the other
indices fixed. The mode-$n$ rank is defined as the dimension of the
subspace spanned by the mode-$n$ vectors of the tensor. The {\em
multilinear rank} of the tensor is the triplet (mode-1 rank, mode-2
rank, mode-3 rank). The mode-$n$ rank generalizes the row and column rank
of matrices. Note that a tensor with multilinear rank $(1,1,1)$ has rank 1
and vice versa. The multilinear rank is invariant under
invertible multilinear transformation \cite[section 2]{DSL}.

Related to the orbits of $2\times 2\times 2$ tensors is the {\em
hyperdeterminant}. Slab operations on $[{\bf X}_1\,|\,{\bf X}_2]$
generate new slabs of the form $\lambda_1\,{\bf
X}_1+\lambda_2\,{\bf X}_2$. There holds
\be
\label{eq-detslabmix}
{\rm det}(\lambda_1\,{\bf X}_1+\lambda_2\,{\bf X}_2)=
\lambda_1^2\,{\det}({\bf X}_1)+\lambda_1\,\lambda_2\,
\frac{{\rm det}({\bf X}_1+{\bf X}_2)-{\rm det}({\bf X}_1-{\bf X}_2)}
{2}\,+\lambda_2^2\,{\rm det}({\bf X}_2)\,.
\ee

\noindent The hyperdeterminant of $\ul{\bf X}$, denoted as $\Delta(\ul{\bf
X})$, is defined as the discriminant of the quadratic polynomial
(\ref{eq-detslabmix}):
\be
\label{eq-hypdet}
\Delta(\ul{\bf X})=\left[
\frac{{\rm det}({\bf X}_1+{\bf X}_2)-{\rm det}({\bf X}_1-{\bf X}_2)}
{2}\right]^2 -4\,{\det}({\bf X}_1)\,{\rm det}({\bf X}_2)\,.
\ee

\noindent Hence, if $\Delta(\ul{\bf X})$ is nonnegative, then a real
slabmix exists that is singular. If $\Delta(\ul{\bf X})$ is positive, then
two real and linearly independent singular slabmixes exist. It follows from
(\ref{eq-detslabmix})-(\ref{eq-hypdet}) that the hyperdeterminant is equal
to the discriminant of the characteristic polynomial of ${\rm det}({\bf
X}_1){\bf X}_2{\bf X}_1^{-1}$ or ${\rm det}({\bf X}_2){\bf X}_1{\bf
X}_2^{-1}$. The sign of the hyperdeterminant is invariant under invertible
multilinear transformation \cite[section 5]{DSL}.

Table~\ref{tab-1} lists the canonical forms for each orbit as well as their
rank, multilinear rank and hyperdeterminant sign. Generic $2\times 2\times
2$ tensors have rank 2 or 3 over the real field, both on a set of positive
measure \cite{Krus89} \cite{TBK}.

For later use, we state the following rank and orbit criteria. The
rank criteria have been proven for $p\times p\times 2$ tensors in
\cite{JJ}. The $2\times 2\times 2$ orbits can be found in \cite[section
7]{DSL}. In the sequel, we will use this result to verify the orbit of a
$2\times 2\times 2$ tensor.

\begin{lem}
\label{lem-DSL}
Let $\ul{\bf X}$ be a $2\times 2\times 2$ tensor with slabs ${\bf X}_1$
and ${\bf X}_2$.
\begin{itemize}
\item[$(i)$] If ${\bf X}_2{\bf X}_1^{-1}$ or ${\bf X}_1{\bf X}_2^{-1}$
has real
eigenvalues and is diagonalizable, then $\ul{\bf X}$ is in orbit $G_2$.
\item[$(ii)$] If ${\bf X}_2{\bf X}_1^{-1}$ or ${\bf X}_1{\bf X}_2^{-1}$
has two identical real eigenvalues with only one associated eigenvector,
then $\ul{\bf X}$ is in orbit $D_3$.
\item[$(iii)$] If ${\bf X}_2{\bf X}_1^{-1}$ or ${\bf X}_1{\bf X}_2^{-1}$
has complex eigenvalues, then $\ul{\bf X}$ is in orbit $G_3$. \end{itemize}
\ep
\end{lem}

\begin{table}[p]
\begin{center}
\begin{tabular}{cccc}
\hline

canonical form & tensor rank & multilinear rank & sign $\Delta$ \\[2mm]

\hline\\

$D_0:\;\left[\begin{array}{cc|cc}0&0&0&0\\0&0&0&0\end{array}\right]$
& 0 & $(0,0,0)$ & 0 \\[7mm]

$D_1:\;\left[\begin{array}{cc|cc}1&0&0&0\\0&0&0&0\end{array}\right]$
& 1 & $(1,1,1)$ & 0 \\[7mm]

$D_2:\;\left[\begin{array}{cc|cc}1&0&0&0\\0&1&0&0\end{array}\right]$
& 2 & $(2,2,1)$ & 0 \\[7mm]

$D'_2:\;\left[\begin{array}{cc|cc}1&0&0&1\\0&0&0&0\end{array}\right]$
& 2 & $(1,2,2)$ & 0 \\[7mm]

$D''_2:\;\left[\begin{array}{cc|cc}1&0&0&0\\0&0&1&0\end{array}\right]$
& 2 & $(2,1,2)$ & 0 \\[7mm]

$G_2:\;\left[\begin{array}{cc|cc}1&0&0&0\\0&0&0&1\end{array}\right]$
& 2 & $(2,2,2)$ & $+$ \\[7mm]

$D_3:\;\left[\begin{array}{cc|cc}0&1&1&0\\1&0&0&0\end{array}\right]$
& 3 & $(2,2,2)$ & 0 \\[7mm]

$G_3:\;\left[\begin{array}{cc|cc}-1&0&0&1\\0&1&1&0\end{array}\right]$
& 3 & $(2,2,2)$ & $-$ \\[7mm]

\hline

\end{tabular}
\end{center}

\caption{Orbits of $2\times 2\times 2$ tensors under the action of
invertible multilinear transformation $({\bf S},{\bf T},{\bf U})$ over the
real field. The letters $D$ and $G$ stand for ``degenerate" (zero volume
set in the 8-dimensional space of $2\times 2\times 2$ tensors) and
``typical" (positive volume set), respectively.} \label{tab-1}
\end{table}

\newpage
\section{Best rank-1 subtraction for $2\times 2\times 2$ tensors}
\setcounter{equation}{0}
For $2\times 2\times 2$ tensors $\ul{\bf X}$ in the orbits of
Table~\ref{tab-1}, we would like to know in which orbit $\ul{\bf X}-\ul{\bf
Y}$ is contained, where $\ul{\bf Y}$ is a best rank-1 approximation of
$\ul{\bf X}$. In this section, we present both examples and general
results. We begin by formulating our main result. It is not a deterministic
result, but involves generic $2\times 2\times 2$ tensors, which are in
orbits $G_2$ and $G_3$. Any tensor randomly generated from a continuous
distribution can be considered to be typical. The full proof of
Theorem~\ref{t-1} is contained in Appendix A.

\begin{thm}
\label{t-1}
Let $\ul{\bf X}$ be a generic $2\times 2\times 2$ tensor, and let $\ul{\bf
Y}$ be a best rank-$1$ approximation of $\ul{\bf X}$. Then almost
all tensors $\ul{\bf X}-\ul{\bf Y}$ are in orbit $D_3$.
\end{thm}

\noindent {\bf Proof sketch.} We proceed as in the first part of Section 2.
We show that there are eight stationary points $({\bf y},{\bf z})$
satisfying (\ref{eq-gradsub2y})-(\ref{eq-gradsub2z}), and that these can be
obtained as roots of an $8$th degree polynomial. There are two stationary
points that yield ${\bf x}={\bf 0}$ in (\ref{eq-xopt}), and do not
correspond to the minimum of criterion (\ref{eq-crit}). For the other six
stationary points, we have $\Delta(\ul{\bf X}-\ul{\bf Y})=0$, where
$\ul{\bf Y}$ is the corresponding rank-1 tensor. Finally, we show that the
multilinear rank of $\ul{\bf X}-\ul{\bf Y}$ equals $(2,2,2)$ for these six
rank-1 tensors $\ul{\bf Y}$. Hence, it follows that the best rank-1
approximation $\ul{\bf Y}$ satisfies $\Delta(\ul{\bf X}-\ul{\bf Y})=0$ and
that the multilinear rank of $\ul{\bf X}-\ul{\bf Y}$ is equal to $(2,2,2)$.
From Table~\ref{tab-1} it then follows that $\ul{\bf X}-\ul{\bf Y}$ is in
orbit $D_3$. \ep

\noindent Hence, for typical tensors in orbit $G_2$, subtracting a best
rank-1 approximation {\em increases} the rank to 3. For typical tensors in
orbit $G_3$, subtracting a best rank-1 approximation does not affect the
rank. This is completely different from matrix analysis.

In the proof of Theorem~\ref{t-1} in Appendix A, it is shown that the slabs
of $\ul{\bf Z}=\ul{\bf X}-\ul{\bf Y}$ are nonsingular almost everywhere.
From Lemma~\ref{lem-DSL} it follows that ${\bf Z}_2{\bf Z}_1^{-1}$ has
identical real eigenvalues and is not diagonalizable, while ${\bf X}_2{\bf
X}_1^{-1}$ has either distinct real eigenvalues or complex eigenvalues.
Hence, the subtraction of a best rank-1 approximation yields identical real
eigenvalues.

Next, we consider $\ul{\bf X}$ in other orbits, and present deterministic
results. We have the following result for the degenerate orbits of ranks 1
and 2.

\begin{propos}
\label{p-1}
Let $\ul{\bf X}$ be a $2\times 2\times 2$ tensor, and let $\ul{\bf Y}$ be a
best rank-$1$ approximation of $\ul{\bf X}$.
\begin{itemize}
\item[$(i)$] If $\ul{\bf X}$ is in orbit $D_1$, then $\ul{\bf X}-\ul{\bf
Y}$ is in orbit $D_0$.
\item[$(ii)$] If $\ul{\bf X}$ is in orbit $D_2$, $D_2'$, or $D_2''$, then
$\ul{\bf X}-\ul{\bf Y}$ is in orbit $D_1$.
\end{itemize}
\end{propos}

\noindent {\bf Proof.} For $\ul{\bf X}$ in orbit
$D_1$ it is obvious that $\ul{\bf Y}=\ul{\bf X}$ is the unique best rank-1
approximation. Then $\ul{\bf X}-\ul{\bf Y}$ is in orbit $D_0$.

Next, let $\ul{\bf X}$ be in orbit $D_2$. Then there exist orthonormal
${\bf S},{\bf T},{\bf U}$ such that
\be
({\bf S},{\bf T},{\bf U})\cdot\ul{\bf X}=
\left[\begin{array}{cc|cc}\lambda&0&0&0\\0&\mu&0&0\end{array}\right]\,,
\ee

\noindent see \cite[proof of lemma 8.2]{DSL}. Subtracting a best rank-1
approximation from this tensor results in $\lambda$ or $\mu$ being set to
zero (whichever has the largest absolute value; for $\lambda=\mu$ there are
two best rank-1 approximations). Hence, the result is a rank-1 tensor. From
Lemma~\ref{p-orth} it follows that the same is true for subtracting a best
rank-1 approximation from $\ul{\bf X}$. For $\ul{\bf X}$ in orbits $D_2'$
and $D_2''$ the proof is analogous. \ep

\noindent For $\ul{\bf X}$ in orbit $G_2$ or $D_3$, the tensor
$\ul{\bf X}-\ul{\bf Y}$ is not restricted to a single orbit. The following
examples illustrate this fact.

\begin{exa}
\label{ex-G2}
{\rm Let
\be
\ul{\bf X} = \left[\begin{array}{cc|cc}1 & 0 & 0 & 0 \\ 0 & 0 & 0
& 1\end{array} \right]\,,
\ee

\noindent which is the canonical tensor of
orbit $G_2$ in Table~\ref{tab-1}. It can be seen that $\ul{\bf X}-\ul{\bf
Y}$ is in $D_1$ (the only nonzero entry of $\ul{\bf Y}$ is either
$Y_{111}$ or $Y_{222}$).

On the other hand, consider
\be
\ul{\bf X}=\left[\begin{array}{cc|cc}0 & 1 & 1 & 0 \\ 1 & 0 & 0 &
2\end{array} \right]\,. \ee

\noindent For this tensor, ${\bf X}_2{\bf X}_1^{-1}$ has two distinct real
eigenvalues. Hence, by Lemma~\ref{lem-DSL}, the tensor is in orbit $G_2$.
It can be shown that $\ul{\bf X}$ has a unique
best rank-1 approximation $\ul{\bf Y}$ and that
\be
\ul{\bf X}-\ul{\bf Y}=
\left[\begin{array}{cc|cc}0 & 1 & 1 & 0 \\ 1 & 0 & 0 & 0\end{array}
\right]\,,
\ee

\noindent which is the canonical tensor of orbit $D_3$ in
Table~\ref{tab-1}.}\ep
\end{exa}

\begin{exa}
\label{ex-D3}
{\rm It follows from Lemma~\ref{lem-DSL} that the following
tensors are in orbit $D_3$:
\be
\left[\begin{array}{cc|cc}2 & 0 & 0 & 1 \\ 0 & 0 & 1 &
0\end{array} \right]\,,\quad\quad
\left[\begin{array}{cc|cc}1 & 0 & 0 & 1 \\ 0 & 0 & 2 &
0\end{array} \right]\,,\quad\quad
\left[\begin{array}{cc|cc}1 & 0 & 0 & 2 \\ 0 & 0 & 1 &
0\end{array} \right]\,.
\ee

\noindent Subtracting the best rank-1 approximation $\ul{\bf Y}$ from these
tensors amounts to replacing the element 2 by zero. Hence, $\ul{\bf
X}-\ul{\bf Y}$ is in orbit $D_2$, $D_2'$, and $D_2''$, respectively.}\ep
\end{exa}

\noindent On the other hand, it can be verified numerically or analytically
that for $\ul{\bf X}$ equal to the canonical tensor of orbit $D_3$ in
Table~\ref{tab-1}, we have $\ul{\bf X}-\ul{\bf Y}$ also in orbit $D_3$.
Moreover, numerical experiments show that for a generic $\ul{\bf X}$ in
orbit $D_3$, we have $\ul{\bf X}-\ul{\bf Y}$ in orbit $D_3$ as well. This
suggests the following.

\begin{conjec}
If $\ul{\bf X}$ is in orbit $D_3$ and $\ul{\bf Y}$ is a best rank-$1$
approximation of $\ul{\bf X}$, then almost all tensors $\ul{\bf X}-\ul{\bf
Y}$ are in $D_3$.\ep \end{conjec}

\noindent The tensor $\ul{\bf X}$ in Example~\ref{ex-KHL} is in orbit $G_3$
by Lemma~\ref{lem-DSL}. It can be shown that any of the infinitely
many best rank-1 approximations of $\ul{\bf X}$ yields $\ul{\bf X}-\ul{\bf
Y}$ in orbit $D_3$ (proof available on request). The example below yields
the same result for another $\ul{\bf X}$ in orbit $G_3$. Numerically and
analytically, we have not found any $\ul{\bf X}$ in orbit $G_3$ for which
$\ul{\bf X}-\ul{\bf Y}$ is not in orbit $D_3$.

\begin{exa}
\label{ex-G3}
{\rm Let
\be
\ul{\bf X}=\left[\begin{array}{cc|cc}1&0&0&-2\\0&1&1&0\end{array}\right]\,.
\ee

\noindent Since ${\bf X}_2{\bf X}_1^{-1}$ has complex eigenvalues, $\ul{\bf
X}$ is in orbit $G_3$ by Lemma~\ref{lem-DSL}. It can be verified that
$\ul{\bf X}$ has a unique best rank-1 approximation such that
\be
\ul{\bf X}-\ul{\bf Y}
=\left[\begin{array}{cc|cc}1&0&0&0\\0&1&1&0\end{array}\right]\,.
\ee

\noindent The latter tensor can be transformed to the canonical form of
orbit $D_3$ by swapping rows within each slab (i.e., by applying a
permutation in the first mode).}\ep \end{exa}

\noindent Our next result concerns tensors with diagonal slabs, i.e.
\be
\label{eq-Xdiag}
\ul{\bf X}=\left[\begin{array}{cc|cc}a&0&e&0\\0&d&0&h\end{array}\right]\,.
\ee

\noindent Then $\ul{\bf X}$ has rank at most 2, since
\be
\ul{\bf X}=\left(\begin{array}{c}1\\
0\end{array}\right) \otimes \left(\begin{array}{c}1\\
0\end{array}\right) \otimes \left(\begin{array}{c}a\\
e\end{array}\right) + \left(\begin{array}{c}0\\
1\end{array}\right) \otimes \left(\begin{array}{c}0\\
1\end{array}\right) \otimes \left(\begin{array}{c}d\\
h\end{array}\right)\,.
\ee

\noindent Also, if $ah\neq de$, then $\ul{\bf X}=({\bf I}_2,{\bf I}_2,{\bf
U})\cdot\ul{\widetilde{\bf X}}$, where $\ul{\widetilde{\bf X}}$ is the
canonical tensor of orbit $G_2$ in Table~\ref{tab-1}, and
\be
{\bf U}=\left[\begin{array}{cc}a& d\\
e&h\end{array}\right]\,.
\ee

\noindent Hence, in this case $\ul{\bf X}$ is in orbit $G_2$.

We show that, for $2\times 2\times 2$ tensors with diagonal slabs, we have
$\ul{\bf X}-\ul{\bf Y}$ in orbit $D_1$. Naturally, the same holds for
$\ul{\bf X}$ that can be transformed to a tensor with diagonal slabs by
orthonormal multilinear transformation (see Lemma~\ref{p-orth}). Note that
tensors with diagonal slabs in orbit $G_2$ form an exception to the
result of Theorem~\ref{t-1}, as does the canonical tensor of orbit $G_2$
(see Example~\ref{ex-G2}). However, Theorem~\ref{t-1} states that these
exceptions form a set of measure zero.

\begin{propos}
\label{p-3}
Let $\ul{\bf X}$ be a $2\times 2\times 2$ tensor with diagonal slabs and
rank $2$, and let $\ul{\bf Y}$ be a best rank-$1$ approximation of
$\ul{\bf X}$. Then $\ul{\bf X}-\ul{\bf Y}$ is in orbit $D_1$.
\end{propos}

\noindent {\bf Proof.} We use the first part of Section 2.
Let $\ul{\bf X}$ be as in (\ref{eq-Xdiag}). First, we assume
$a^2+e^2<d^2+h^2$. For
\be
\label{eq-Yopt}
\ul{\bf
Y}=\left[\begin{array}{cc|cc}0&0&0&0\\0&d&0&h\end{array}\right]\,, \ee

\noindent we have $\|\ul{\bf X}-\ul{\bf Y}\|^2=a^2+e^2$ and $\ul{\bf
X}-\ul{\bf Y}$ in orbit $D_1$. Next, we show that (\ref{eq-Yopt}) is the
unique best rank-1 approximation of $\ul{\bf X}$. Using (\ref{eq-critsub}),
the equation $\|\ul{\bf X}-\ul{\bf Y}\|^2\le a^2+e^2$ can be written as
\be
(d^2+h^2)\,(y_1^2+y_2^2)\,(z_1^2+z_2^2)\le
(ay_1z_1+ey_1z_2)^2+(dy_2z_1+hy_2z_2)^2\,,
\ee

\noindent which can be rewritten as
\be
\label{eq-proof1}
(d^2+h^2-a^2-e^2)\,y_1^2\,(z_1^2+z_2^2)+(ez_1-az_2)^2\,y_1^2+
(hz_1-dz_2)^2\,y_2^2\le 0\,.
\ee

\noindent Since $(d^2+h^2-a^2-e^2)$ is positive by assumption, and ${\bf
y}$ nor ${\bf z}$ can be all-zero, it follows that (\ref{eq-proof1}) can
only hold with equality, that is, for $y_1=0$ and $hz_1=dz_2$. Using
(\ref{eq-xopt}), it then follows that the $\ul{\bf Y}$ for which we have
equality in (\ref{eq-proof1}) is given by (\ref{eq-Yopt}). This shows that
(\ref{eq-Yopt}) is the unique best rank-1 approximation of $\ul{\bf X}$.

Next, we consider the case $a^2+e^2>d^2+h^2$. Analogous to the first part
of the proof, it can be shown that
\be
\label{eq-Yopt2}
\ul{\bf
Y}=\left[\begin{array}{cc|cc}a&0&e&0\\0&0&0&0\end{array}\right]\,, \ee

\noindent is the unique best rank-1 approximation of $\ul{\bf X}$.
This implies that $\|\ul{\bf X}-\ul{\bf Y}\|^2=d^2+h^2$ and $\ul{\bf
X}-\ul{\bf Y}$ is in orbit $D_1$.

Finally, we consider the case $a^2+e^2=d^2+h^2$. Here, we have multiple
best rank-1 approximations. Setting $y_1=0$ in (\ref{eq-proof1}) yields
(\ref{eq-Yopt}) as a best rank-1 approximation. Setting $y_2=0$ yields
(\ref{eq-Yopt2}) as a best rank-1 approximation.
If $ah=de$, then (\ref{eq-proof1}) can also be satisfied by
setting $ez_1=az_2$ and $hz_1=dz_2$. This yields
\be
\label{eq-Yopt3}
\ul{\bf Y}=(y_1^2+y_2^2)^{-1}\;
\left[\begin{array}{cc|cc}y_1^2a&y_1y_2a&y_1^2e&y_1y_2e\\y_1y_2d&y_2^2d&y_1y
_2h&y_2^2h\end{array}\right]\,. \ee

\noindent It can be verified that for $\ul{\bf Y}$ in (\ref{eq-Yopt3}) we
have $\|\ul{\bf X}-\ul{\bf Y}\|^2=d^2+h^2=a^2+e^2$ and $\ul{\bf
X}-\ul{\bf Y}$ in orbit $D_1$. This completes the proof.
\ep

\section{Best rank-1 approximation for symmetric tensors}
\setcounter{equation}{0}
Here, we consider the best rank-1 approximation problem for a $3$rd
order tensor $\ul{\bf X}\in\R^{d\times d\times d}$ that is symmetric
in all modes, i.e. $X_{ijk}=X_{jik}=X_{kji}=X_{ikj}=X_{jki}=X_{kij}$. We
assume the same for the rank-1 approximation, which yields the problem \be
\label{eq-prob2} {\rm min}_{{\bf y}\in\R^{d}}\;
||\ul{\bf X}-{\bf y}\otimes{\bf y}\otimes{\bf y}||^2\,. \ee

\noindent An adaption of \cite[proposition 4.2]{DSL} shows that problem
(\ref{eq-prob2}) always has an optimal solution.

Analogous to the first part of Section 2, the criterion of (\ref{eq-prob2})
can be written as
\be
\label{eq-crit2}
\Psi_2=||\ul{\bf X}||^2 -2\,\ul{\bf X}\bullet_1 {\bf y}^T\bullet_2 {\bf y}^T
\bullet_3 {\bf y}^T + ||{\bf y}||^6 \,.
\ee

\noindent When the gradient of $\Psi_2$ with respect to ${\bf y}$ is set
to zero, we obtain
\be
\label{eq-grad2} {\bf y}  =
\frac{\ul{\bf X} \bullet_1 {\bf y}^T \bullet_2 {\bf y}^T}{||{\bf y}||^4}\,.
\ee

\noindent Substituting (\ref{eq-grad2}) into the
criterion (\ref{eq-crit2}) yields
\be
\label{eq-critsub2}
\Psi_2 = ||\ul{\bf X}||^2 -\frac{||\ul{\bf X}\bullet_1 {\bf y}^T\bullet_2
{\bf y}^T||^2} {||{\bf y}||^4}\,. \ee

\noindent Hence, a best symmetric rank-1 approximation $({\bf y}\otimes {\bf
y}\otimes {\bf y})$ of $\ul{\bf X}$ is found by minimizing (\ref{eq-crit2})
or (\ref{eq-critsub2}) over ${\bf y}$. The stationary points ${\bf y}$ are
given by (\ref{eq-grad2}); this was already noticed in \cite[section
2.3]{Com98}.

\section{Rank criteria and orbits of symmetric $2\times 2\times 2$ tensors}
\setcounter{equation}{0}
Here, we consider real symmetric $2\times 2\times 2$ tensors. We establish
their ranks and orbits under invertible multilinear transformation $({\bf
S},{\bf S},{\bf S})\cdot\ul{\bf X}$. These transformations preserve the
symmetry. The {\em symmetric tensor rank} \cite[section 4]{CGLM} is defined
as the smallest $R$ such that
\be
\label{eq-symdecomp}
\ul{\bf X}=\sum_{r=1}^R {\bf a}_r \otimes
{\bf a}_r \otimes {\bf a}_r\,. \ee

\noindent There is a bijection between symmetric $d\times d\times d$
tensors and homogeneous polynomials of degree 3 in $d$ variables. A
symmetric $d\times d\times d$ tensor $\ul{\bf X}$ is associated with the
polynomial
\be
\label{eq-pol}
p(u_1,\ldots,u_d)=\sum_{ijk} x_{ijk}\;u_i u_j u_k\,.
\ee

\noindent A multilinear transformation $({\bf S},{\bf S},{\bf
S})\cdot\ul{\bf X}$ is equivalent to a change of variables ${\bf v}={\bf
Su}$ in the associated polynomial.

The symmetric rank of symmetric $2\times 2\times\cdots\times 2$
tensors can be obtained from the well-known Sylvester Theorem, which makes
use of the polynomial representation \cite[section 5]{KR84} \cite{ComMou}.
For generic symmetric $2\times 2\times\cdots\times 2$ tensors,
\cite{ComMou} show that the Sylvester Theorem defines an algorithm to
compute a symmetric decomposition (\ref{eq-symdecomp}) with $R$ equal to
the symmetric rank. Below, the Sylvester Theorem for symmetric $2\times
2\times 2$ tensors is formulated.

\begin{thm}
\label{t-sylv}
{\rm\bf (Sylvester)} A real symmetric $2\times 2\times 2$ tensor with
associated polynomial
\be
p(u_1,u_2)=\gamma_3\,u_1^3+3\,\gamma_2\,u_1^2u_2+3\,\gamma_1\,
u_1u_2^2+\gamma_0\,u_2^3\,,
\ee

\noindent has a symmetric decomposition $(\ref{eq-symdecomp})$ into $R$
rank-1 terms if and only if there exists a vector ${\bf
g}=(g_0,\ldots,g_R)^T$ with
\be
\label{eq-sylv}
\left[\begin{array}{ccc}
\gamma_0 & \ldots & \gamma_R \\[2mm]
\gamma_1 & \ldots & \gamma_{R+1} \\[2mm]
\vdots & & \vdots \\[2mm]
\gamma_{3-R} & \ldots & \gamma_3 \end{array}\right]\;{\bf g}={\bf 0}\,,
\ee

\noindent and if the polynomial
$q(u_1,u_2)=g_R\,u_1^R+\cdots+g_1\,u_1u_2^{R-1}+g_0\,u_2^R$ has $R$
distinct real roots. \ep \end{thm}

\noindent For our purposes, we make use of a symmetric rank criterion
similar to Lemma~\ref{lem-DSL}, formulated as Lemma~\ref{lem-symrank}
below. The link between this rank criterion and the Sylvester Theorem will
be explained at the end of this section.

Let the entries of a symmetric $2\times 2\times 2$ tensor be denoted as
\be
\label{eq-Xgsym}
\ul{\bf
X}=\left[\begin{array}{cc|cc}a&b&b&c\\b&c&c&d\end{array}\right]\,. \ee

\noindent For later use, we mention that the hyperdeterminant
(\ref{eq-hypdet}) of $\ul{\bf X}$ in (\ref{eq-Xgsym}) is given by
\be
\label{eq-hypdetsym}
\Delta(\ul{\bf X})=(bc-ad)^2-4(bd-c^2)(ac-b^2)\,.
\ee

\noindent As in the asymmetric case, the sign of the hyperdeterminant is
invariant under invertible multilinear transformation $({\bf
S},{\bf S},{\bf S})\cdot\ul{\bf X}$.

\begin{lem}
\label{lem-symrank}
Let $\ul{\bf X}$ be a real symmetric $2\times 2\times 2$ tensor with
slabs ${\bf X}_1$ and ${\bf X}_2$. \begin{itemize} \item[$(i)$]
If ${\bf X}_2{\bf X}_1^{-1}$ or ${\bf X}_1{\bf X}_2^{-1}$ has distinct real
eigenvalues, then $\ul{\bf X}$ has symmetric rank $2$.
\item[$(ii)$] If ${\bf X}_2{\bf X}_1^{-1}$ or ${\bf X}_1{\bf
X}_2^{-1}$ has identical real eigenvalues, then $\ul{\bf X}$ has symmetric
rank at least $3$. \item[$(iii)$] If ${\bf X}_2{\bf X}_1^{-1}$ or ${\bf
X}_1{\bf X}_2^{-1}$ has complex eigenvalues, then $\ul{\bf X}$ has
symmetric rank at least $3$. \end{itemize} \end{lem}

\noindent {\bf Proof.} First, we prove $(i)$. We consider ${\bf X}_2{\bf
X}_1^{-1}$. The proof for ${\bf X}_1{\bf X}_2^{-1}$ is completely
analogous. Note that since ${\bf X}_1$ is nonsingular, the symmetric rank
of $\ul{\bf X}$ is at least 2. Let ${\bf X}_2{\bf X}_1^{-1}$ have distinct
real eigenvalues $\lambda_1$ and $\lambda_2$. Using (\ref{eq-Xgsym}), we
have
\be
\label{eq-X2X1}
{\bf X}_2{\bf X}_1^{-1}=\left[\begin{array}{cc}
0 & 1 \\
x & y\end{array}\right]\,,\quad\quad
x=\frac{c^2-bd}{ac-b^2}\,,\quad\quad
y=\frac{ad-bc}{ac-b^2}\,,
\ee

\noindent which has eigenvalues
\be
\label{eq-eig12}
\lambda_1=\frac{y+\sqrt{y^2+4x}}{2}\,,\quad\quad\quad\quad
\lambda_2=\frac{y-\sqrt{y^2+4x}}{2}\,.
\ee

\noindent It can be verified that the eigenvectors of ${\bf X}_2{\bf
X}_1^{-1}$ are $\alpha\,(1\;\;\lambda_1)^T$ and
$\beta\,(1\;\;\lambda_2)^T$, respectively. Next, we show that appropriate
choices of $\alpha$ and $\beta$ yield a symmetric rank-2 decomposition
(\ref{eq-symdecomp}) for $\ul{\bf X}$.

Let ${\bf A}$ contain the eigenvectors of ${\bf X}_2{\bf X}_1^{-1}$, i.e.
\be
\label{eq-A}
{\bf A}=\left[\begin{array}{cc}
\alpha & \beta \\
\alpha\lambda_1 & \beta\lambda_2\end{array}\right]\,.
\ee

\noindent Note that a symmetric rank-2 decomposition (\ref{eq-symdecomp})
can be denoted as ${\bf X}_k={\bf A}\,{\bf A}_k\,{\bf A}^T$, with ${\bf
A}_k={\rm diag}(a_{k1},a_{k2})$, $k=1,2$. The eigendecomposition of
${\bf X}_2{\bf X}_1^{-1}$ is then given by ${\bf A}\,{\bf A}_2{\bf
A}_1^{-1}\,{\bf A}^{-1}$, which is consistent with (\ref{eq-A}) since
${\bf A}_2{\bf A}_1^{-1}={\rm diag}(\lambda_1,\lambda_2)$. To obtain a
symmetric rank-2 decomposition, it remains to solve the equations
\be
\label{eq-remainsym}
{\bf A}^{-1}\,{\bf X}_1={\bf A}_1\,{\bf A}^T=\left[\begin{array}{cc}
\alpha^3 & \alpha^3\lambda_1\\
\beta^3 & \beta^3\lambda_2\end{array}\right]\,,\quad\quad
{\bf A}^{-1}\,{\bf X}_2={\bf A}_2\,{\bf A}^T=\left[\begin{array}{cc}
\alpha^3\lambda_1 & \alpha^3\lambda_1^2\\
\beta^3\lambda_2 & \beta^3\lambda_2^2\end{array}\right]\,.
\ee

\noindent By writing out the entries of ${\bf A}^{-1}{\bf X}_1$
and ${\bf A}^{-1}{\bf X}_2$, it can be seen that a solution
$(\alpha,\beta)$ exists if and only if \be
\label{eq-needthisnow}
a\,\lambda_1\,\lambda_2-b\,(\lambda_1+\lambda_2)+c=0\,,\quad\quad\quad
b\,\lambda_1\,\lambda_2-c\,(\lambda_1+\lambda_2)+d=0\,. \ee

\noindent Using (\ref{eq-eig12}), we have $\lambda_1\lambda_2=-x$ and
$\lambda_1+\lambda_2=y$. Combined with the expressions for $x$ and
$y$ in (\ref{eq-X2X1}), this verifies (\ref{eq-needthisnow}). This completes
the proof of $(i)$.

As shown in the proof of $(i)$, if $\ul{\bf X}$ has one slab nonsingular,
say ${\bf X}_1$, and symmetric rank 2, then ${\bf X}_2{\bf X}_1^{-1}$ is
diagonalizable. As also shown in the proof of $(i)$, when ${\bf X}_2{\bf
X}_1^{-1}$ has identical real eigenvalues it does not have two linearly
independent eigenvectors and, hence, is not diagonalizable. Therefore, in
case $(ii)$ the symmetric rank of $\ul{\bf X}$ is at least 3. The same
holds in case $(iii)$. This completes the proof.\ep

\begin{table}[tp]
\begin{center}
\begin{tabular}{cccc}
\hline

canonical form & symmetric rank & multilinear rank & sign $\Delta$ \\[2mm]

\hline\\

$D_0:\;\left[\begin{array}{cc|cc}0&0&0&0\\0&0&0&0\end{array}\right]$
& 0 & $(0,0,0)$ & 0 \\[7mm]

$D_1:\;\left[\begin{array}{cc|cc}1&0&0&0\\0&0&0&0\end{array}\right]$
& 1 & $(1,1,1)$ & 0 \\[7mm]

$G_2:\;\left[\begin{array}{cc|cc}1&0&0&0\\0&0&0&1\end{array}\right]$
& 2 & $(2,2,2)$ & $+$ \\[7mm]

$D_3:\;\left[\begin{array}{cc|cc}0&1&1&0\\
1 & 0 & 0 & 0\end{array}\right]$
& 3 & $(2,2,2)$ & 0 \\[7mm]

$G_3:\;\left[\begin{array}{cc|cc}-1&0&0&1\\
0 & 1 & 1 & 0\end{array}\right]$
& 3 & $(2,2,2)$ & $-$ \\[7mm]

\hline

\end{tabular}
\end{center}

\caption{Orbits of symmetric $2\times 2\times 2$ tensors under the action of
invertible multilinear transformation $({\bf S},{\bf S},{\bf S})$ over the
real field. The letters $D$ and $G$ stand for ``degenerate" (zero volume
set in the 4-dimensional space of symmetric $2\times 2\times 2$ tensors)
and ``typical" (positive volume set), respectively.} \label{tab-orbit2}
\end{table}

\begin{propos}
\label{p-orbit2}
The orbits of real symmetric $2\times 2\times 2$ tensors $\ul{\bf X}$
under the action of invertible multilinear
transformation $({\bf S},{\bf S},{\bf S})\cdot\ul{\bf X}$, are as given in
Table~$\ref{tab-orbit2}$. \end{propos}

\noindent {\bf Proof.} Let ${\bf e}_1$ and ${\bf e}_2$ denote the first and
second column of ${\bf I}_2$. Orbit $D_0$ corresponds to the all-zero
tensor. For $\ul{\bf X}$ with symmetric rank 1, we have $\ul{\bf X}={\bf
a}\otimes{\bf a}\otimes{\bf a}$. There exists a nonsingular ${\bf S}$ with
${\bf Sa}={\bf e}_1$. Then $({\bf S},{\bf S},{\bf S})\cdot\ul{\bf X}$
equals the canonical form of orbit $D_1$ in Table~\ref{tab-orbit2}, which
has multilinear rank $(1,1,1)$. Analogously, for $\ul{\bf X}$ with
symmetric rank 2, we have $\ul{\bf X}={\bf a}_1\otimes{\bf a}_1\otimes{\bf
a}_1+{\bf a}_2\otimes{\bf a}_2\otimes{\bf a}_2$, with ${\bf a}_1$ and ${\bf
a}_2$ linearly independent. There exists a nonsingular ${\bf S}$ with ${\bf
S}[{\bf a}_1\;{\bf a}_2]=[{\bf e}_1\;{\bf e}_2]$. Then $({\bf S},{\bf
S},{\bf S})\cdot\ul{\bf X}$ equals the canonical form of orbit $G_2$ in
Table~\ref{tab-orbit2}, which has positive hyperdeterminant
(\ref{eq-hypdetsym}) and multilinear rank $(2,2,2)$.

Next, let $\ul{\bf X}$ have symmetric rank 3 and decomposition
(\ref{eq-symdecomp}) with ${\bf A}=[{\bf a}_1\;{\bf a}_2\;{\bf a}_3]$. No
two columns of ${\bf A}$ are proportional, since otherwise a symmetric
rank-2 decomposition is possible. It follows that there exists a
nonsingular ${\bf S}$ with
\be
\label{eq-rank3decomp}
{\bf SA}=\left[\begin{array}{ccc}
\alpha & 0 & 1\\
0 & \beta & 1\end{array}\right]\,.
\ee

\noindent This yields
\be
\label{eq-oftheform}
({\bf S},{\bf S},{\bf S})\cdot\ul{\bf X}=\left[\begin{array}{cc|cc}
a & 1 & 1 & 1\\
1 & 1 & 1 & d\end{array}\right]\,,
\ee

\noindent with $a=1+\alpha^3$ and $d=1+\beta^3$. We define orbits $D_3$ and
$G_3$ according to whether the hyperdeterminant $\Delta$ is zero or
negative, respectively. Note that $\Delta>0$ is associated with orbit $G_2$.

In Appendix C, we show that any pair of tensor $\ul{\bf X}$ in orbit $D_3$
or $G_3$ is related to the canonical form $\ul{\bf Y}$ of the orbit by an
invertible multilinear transformation $({\bf S},{\bf S},{\bf
S})\cdot\ul{\bf Y}=\ul{\bf X}$.

We conclude our proof by showing that the symmetric rank of real symmetric
$2\times 2\times 2$ tensors is at most 3. Let $\ul{\bf X}$ be as in
(\ref{eq-Xgsym}). Suppose $b\neq 0$ and $c\neq 0$. Then $({\bf S},{\bf
S},{\bf S})\cdot\ul{\bf X}$ is of the form (\ref{eq-oftheform}) for ${\bf
S}={\rm diag}(\mu,\eta)$ with $\mu^3=c/b^2$ and $\eta^3=b/c^2$.
Since (\ref{eq-oftheform}) has the symmetric rank-3 decomposition
(\ref{eq-rank3decomp}), the tensor $\ul{\bf X}$ has at most symmetric rank
3.

Next, suppose $b=0$ and $c\neq 0$. We subtract a symmetric rank-1 tensor
${\bf a}\otimes{\bf a}\otimes{\bf a}$ with ${\bf a}=\gamma\,{\bf e}_1$ from
$\ul{\bf X}$ such that $(a-\gamma^3)c>0$. Denote the resulting tensor by
$\ul{\bf Z}$. It can be verified that ${\bf Z}_2{\bf Z}_1^{-1}$ has
distinct real eigenvalues. By Lemma~\ref{lem-symrank} $(i)$ it has
symmetric rank 2. Combined with the subtracted rank-1 tensor, this implies
a symmetric rank-3 decomposition of $\ul{\bf X}$.

The case $c=0$ and $b\neq 0$ can be dealt with analogously. When $b=c=0$, a
symmetric rank-2 decomposition is immediate. Hence, the symmetric rank is at
most 3.\ep

\noindent The following corollary follows from Lemma~\ref{lem-symrank} and
Proposition~\ref{p-orbit2}. It is the full analogue of Lemma~\ref{lem-DSL}.

\begin{cor}
\label{cor-symrank}
Let $\ul{\bf X}$ be a real symmetric $2\times 2\times 2$ tensor with
slabs ${\bf X}_1$ and ${\bf X}_2$. \begin{itemize} \item[$(i)$]
If ${\bf X}_2{\bf X}_1^{-1}$ or ${\bf X}_1{\bf X}_2^{-1}$ has distinct real
eigenvalues, then $\ul{\bf X}$ is in orbit $G_2$.
\item[$(ii)$] If ${\bf X}_2{\bf X}_1^{-1}$ or ${\bf X}_1{\bf
X}_2^{-1}$ has identical real eigenvalues, then $\ul{\bf X}$ is in orbit
$D_3$. \item[$(iii)$] If ${\bf X}_2{\bf X}_1^{-1}$ or ${\bf
X}_1{\bf X}_2^{-1}$ has complex eigenvalues, then $\ul{\bf X}$ is in orbit
$G_3$. \end{itemize} \end{cor}

\noindent {\bf Proof.} Since there is only one orbit with symmetric rank
2, the proof of $(i)$ is the proof of Lemma~\ref{lem-symrank} $(i)$.
Since the symmetric rank is at most 3, cases $(ii)$ and $(iii)$ have
symmetric rank 3. As in the asymmetric case, the hyperdeterminant
(\ref{eq-hypdetsym}) is equal to the discriminant of the characteristic
polynomial of ${\rm det}({\bf X}_1){\bf X}_2{\bf X}_1^{-1}$ or ${\rm
det}({\bf X}_2){\bf X}_1{\bf X}_2^{-1}$. Hence, case $(ii)$ has $\Delta=0$
and corresponds to orbit $D_3$, and case $(iii)$ has $\Delta<0$ and
corresponds to orbit $G_3$. \ep

\noindent The hyperdeterminant (\ref{eq-hypdetsym}) is equal to the
discriminant of the polynomial $q(u_1,u_2)$ in the Sylvester Theorem
(Theorem~\ref{t-sylv}) for $\ul{\bf X}$ in (\ref{eq-Xgsym}) and $R=2$.
Indeed, we have $\gamma_3=a$, $\gamma_2=b$, $\gamma_1=c$ and $\gamma_0=d$.
The vector ${\bf g}=(g_0,g_1,g_2)^T$ should satisfy
\be
\left[\begin{array}{ccc}
d & c & b \\
c & b & a \end{array}\right]\;{\bf g}={\bf 0}\,,
\ee

\noindent which implies
\be
g_0=ac-b^2\,,\quad\quad\quad
g_1=bc-ad\,,\quad\quad\quad
g_2=bd-c^2\,.
\ee

\noindent The discriminant of $q(u_1,u_2)$ is given by $g_1^2-4\,g_0\,g_2$
which is equal to the hyperdeterminant (\ref{eq-hypdetsym}). This
establishes the equivalence between the Sylvester Theorem with $R=2$ and
the symmetric rank criteria of Lemma~\ref{lem-symrank}.

\section{Best rank-1 subtraction for symmetric $2\times 2\times 2$ tensors}
\setcounter{equation}{0}
Here, we consider the problem of determining the rank and orbit of $\ul{\bf
X}-\ul{\bf Y}$, where $\ul{\bf X}$ is a symmetric $2\times 2\times 2$
tensor and ${\bf Y}$ is a best symmetric rank-1 approximation of $\ul{\bf
X}$. Obviously, if $\ul{\bf X}$ is in orbit $D_1$, then $\ul{\bf X}-\ul{\bf
Y}$ is in orbit $D_0$. Next, we present our main result in
Theorem~\ref{t-2}, which is the analogue of Theorem~\ref{t-1}. It
concerns generic symmetric $2\times 2\times 2$ tensors, which are in
orbits $G_2$ and $G_3$. The full proof of Theorem~\ref{t-2} is contained in
Appendix B.

\begin{thm}
\label{t-2}
Let $\ul{\bf X}$ be a generic symmetric $2\times 2\times 2$ tensor, and let
$\ul{\bf Y}$ be a best symmetric rank-$1$ approximation of $\ul{\bf X}$.
Then almost all tensors $\ul{\bf X}-\ul{\bf Y}$ are in orbit $D_3$. \end{thm}

\noindent {\bf Proof sketch.} We proceed as in Section 5. It is shown in
\cite[section 3.5]{LMVrank1} that there are three stationary points ${\bf
y}$ satisfying (\ref{eq-grad2}), and that these can be obtained as roots of
a $3$rd degree polynomial. We show that, for all three stationary points,
we have $\Delta(\ul{\bf X}-\ul{\bf Y})=0$, where $\ul{\bf Y}$ is the
corresponding rank-1 tensor. Finally, we show that the multilinear rank of
$\ul{\bf X}-\ul{\bf Y}$ equals $(2,2,2)$ for these three rank-1 tensors
$\ul{\bf Y}$. This suffices to conclude that $\ul{\bf X}-\ul{\bf Y}$ is in
orbit $D_3$. \ep

\noindent Hence, as in the asymmetric $2\times 2\times 2$ case, for typical
symmetric tensors in orbit $G_2$, subtracting a best symmetric rank-1
approximation {\em increases} the symmetric rank to 3. For typical
symmetric tensors in orbit $G_3$, subtracting a best symmetric rank-1
approximation does not affect the symmetric rank.

In the proof of Theorem~\ref{t-2} in Appendix B, it is shown that the slabs
of $\ul{\bf Z}=\ul{\bf X}-\ul{\bf Y}$ are nonsingular almost everywhere.
From Lemma~\ref{lem-symrank} it follows that ${\bf Z}_2{\bf Z}_1^{-1}$ has
identical real eigenvalues, while ${\bf X}_2{\bf X}_1^{-1}$ has either
distinct real eigenvalues or complex eigenvalues. Hence, also for
symmetric $2\times 2\times 2$ tensors, the subtraction of a best rank-1
approximation yields identical real eigenvalues.

We conclude this section with examples of $\ul{\bf X}$ in
orbits $G_2$ and $G_3$ such that $\ul{\bf X}-\ul{\bf Y}$ is in orbit $D_3$.

\begin{exa}
\label{ex-sym1}
{\rm Let
\be
\label{eq-Xsym}
\ul{\bf
X}=\left[\begin{array}{cc|cc}0&1&1&1\\1&1&1&0\end{array}\right]\,,
\quad\quad\quad
{\bf X}_2{\bf X}_1^{-1}=\left[\begin{array}{cc}
0&1\\-1&1\end{array}\right]\,.
\end{equation}

\noindent Since the latter has complex eigenvalues, Lemma~\ref{lem-symrank}
shows that $\ul{\bf X}$ is in orbit $G_3$.

Next, we compute the best symmetric rank-1 approximation $\ul{\bf Y}$ to
$\ul{\bf X}$, which has the form
\be
\label{eq-Ysym}
\ul{\bf
Y}=\left[\begin{array}{cc|cc}y_1^3&y_1^2y_2&y_1^2y_2&y_1y_2^2\\y_1^2y_2 &
y_1y_2^2&y_1y_2^2&y_2^3 \end{array}\right]\,. \ee

\noindent The stationary points (\ref{eq-grad2}) are given by
\begin{eqnarray}
\label{eq-der1}
6\,y_1^5+12\,y_1^3y_2^2+6\,y_1y_2^4-12\,y_1y_2 -6\,y_2^2 &=& 0\,,\\[2mm]
\label{eq-der2}
6\,y_2^5+12\,y_1^2y_2^3+6\,y_1^4y_2-12\,y_1y_2 -6\,y_1^2 &=& 0\,.
\end{eqnarray}

\noindent It follows that $y_1\neq 0$ and $y_2\neq 0$ (if one of them
equals zero, then both are zero and $\ul{\bf Y}$ is all-zero). Muliplying
(\ref{eq-der2}) by $y_1$ and subtracting (\ref{eq-der1}) multiplied by
$y_2$ yields
\be
6\;\large(y_2-y_1\large)\left(y_2-\frac{-3+\sqrt{5}}{2}\,y_1\right)
\left(y_2-\frac{-3-\sqrt{5}}{2}\,y_1\right)=0\,.
\ee

\noindent Hence, $y_2=y_1$ or $y_2=(-3/2\pm\sqrt{5}/2)\,y_1$. However,
it can be verified that the latter is in contradiction with
(\ref{eq-der1})-(\ref{eq-der2}). Therefore, $y_2=y_1$ and it follows from
(\ref{eq-der1})-(\ref{eq-der2}) that $y_1^3=y_2^3=3/4$. This yields
$\Psi=3/2$ and \be \label{eq-Ysymsol} \ul{\bf
Y}=\frac{1}{4}\;\left[\begin{array}{cc|cc}3&3&3&3\\3& 3&3&3
\end{array}\right]\,,\quad\quad\quad
\ul{\bf Z}=
\ul{\bf X}-\ul{\bf Y}
=\frac{1}{4}\;\left[\begin{array}{cc|cc}-3&1&1&1\\1&1&1&-3
\end{array}\right]\,.\ee

\noindent Hence,
\be
{\bf Z}_2{\bf Z}_1^{-1}=\left[\begin{array}{cc}
0&1\\-1&-2\end{array}\right]\,,
\ee

\noindent which has a double eigenvalue $-1$. Hence, $\ul{\bf Z}$ is in
orbit $D_3$ by Lemma~\ref{lem-symrank}.}\ep \end{exa}

\begin{exa}
\label{ex-sym2}
{\rm Let
\be
\label{eq-Xsym2}
\ul{\bf X}=\left[\begin{array}{cc|cc}3&1&1&1\\1&1&1&3\end{array}\right]\,,
\quad\quad\quad
{\bf X}_2{\bf X}_1^{-1}=\left[\begin{array}{cc}
0&1\\-1&4\end{array}\right]\,.
\ee

\noindent Since the latter has real and distinct eigenvalues,
Lemma~\ref{lem-symrank} shows that $\ul{\bf x}$ is in orbit $G_2$.

Analogous to Example~\ref{ex-sym1}, it can be shown that the best
symmetric rank-1 approximation of $\ul{\bf X}$ is given by
\be
\ul{\bf Y}=\frac{3}{2}\;\left[\begin{array}{cc|cc}1&1&1&1\\1&1&1&1
\end{array}\right]\,.
\ee

\noindent We obtain
\be
\ul{\bf Z}=
\ul{\bf X}-\ul{\bf Y}
=\frac{1}{2}\;\left[\begin{array}{cc|cc}3&-1&-1&-1\\-1&-1&-1&3
\end{array}\right]\,,\quad\quad\quad
{\bf Z}_2{\bf Z}_1^{-1}=\left[\begin{array}{cc}
0&1\\-1&-2\end{array}\right]\,.
\ee

\noindent The latter has a double eigenvalue $-1$. Hence,
$\ul{\bf Z}$ is in orbit $D_3$ by Lemma~\ref{lem-symrank}.}\ep
\end{exa}

\section{Discussion}
\setcounter{equation}{0}
It is now rather well known that consecutively subtracting a best
rank-1 approximation from a higher-order tensor generally does not either
reveal tensor rank nor yield a ``good'' low-rank approximation. A
numerical example and discussion is provided in \cite[section 7]{KofReg}.
Hence, a rank-1 deflation procedure as is available for matrices, generally
does not exist for higher-order tensors. We have given a mathematical
treatment of this property for real $2\times 2\times 2$ tensors. In
Theorem~\ref{t-1}, we showed that subtracting a best rank-1 approximation
from a generic $2\times 2\times 2$ tensor (which has rank 2 or 3) results
in a rank-3 tensor located on the boundary between the sets of rank-2
and rank-3 tensors. Hence, for typical tensors of rank 2, subtracting a
best rank-1 approximation {\em increases} the rank to 3.

A generic $2\times 2\times 2$ tensor $\ul{\bf X}$ has rank 2 if ${\bf
X}_2{\bf X}_1^{-1}$ has distinct real eigenvalues, and rank 3 if ${\bf
X}_2{\bf X}_1^{-1}$ has complex eigenvalues; see Lemma~\ref{lem-DSL}. If
$\ul{\bf Y}$ is a best rank-1 approximation of $\ul{\bf X}$, then $\ul{\bf
Z}=\ul{\bf X}-\ul{\bf Y}$ has rank 3 and lies on the boundary between the
rank-2 and rank-3 sets, i.e. ${\bf Z}_2{\bf Z}_1^{-1}$ has identical real
eigenvalues. The rank-2 and rank-3 orbits $G_2$ and $G_3$ are characterized
by positive and negative hyperdeterminant $\Delta$, respectively, while on
the boundary we have $\Delta=0$. The result that subtraction of a best
rank-1 approximation yields identical real eigenvalues for ${\bf Z}_2{\bf
Z}_1^{-1}$ is new and expands the knowledge of the topology of tensor rank.

Numerical experiments yield the conjecture that for a generic real-valued
$p\times p\times 2$ tensor $\ul{\bf X}$, subtracting its best rank-1
approximation $\ul{\bf Y}$ results in $\ul{\bf Z}=\ul{\bf X}-\ul{\bf Y}$
with ${\bf Z}_2{\bf Z}_1^{-1}$ having one pair of identical real
eigenvalues with only one associated eigenvector. Moreover, if the number of
pairs of complex eigenvalues of ${\bf X}_2{\bf X}_1^{-1}$ equals $n$, then
${\bf Z}_2{\bf Z}_1^{-1}$ has $\max(0,n-1)$ pairs of complex eigenvalues.
For $n=0$, this implies that $\ul{\bf X}$ has rank $p$ and $\ul{\bf Z}$ has
rank $p+1$ \cite{JJ} \cite[lemma 2.2]{SDL}.

We also considered real symmetric $2\times 2\times 2$ tensors. In
Lemma~\ref{lem-symrank}, we provided a symmetric rank criterion via the
eigenvalues of ${\bf X}_2{\bf X}_1^{-1}$, which is similar to the asymmetric case.
Symmetric tensors have rank 2 and 3 on sets of positive volume, and ${\bf
X}_2{\bf X}_1^{-1}$ with distinct real eigenvalues implies $\Delta>0$ and
orbit $G_2$, while ${\bf X}_2{\bf X}_1^{-1}$ with complex eigenvalues
implies $\Delta<0$ and orbit $G_3$. When ${\bf X}_2{\bf X}_1^{-1}$ has
identical real eigenvalues, it has $\Delta=0$ and symmetric rank 3 (orbit
$D_3$). The rank criteria of Lemma~\ref{lem-symrank} are equivalent to
the well-known Sylvester Theorem for symmetric rank 2. In Theorem~\ref{t-2},
we showed that subtracting a best symmetric rank-1 approximation from a
typical symmetric tensor yields a tensor in orbit $D_3$, i.e. it has
symmetric rank 3 and $\Delta=0$. This result is completely analogous to the
asymmetric $2\times 2\times 2$ case.

A third case not reported here is that of $2\times 2\times 2$ tensors with
symmetric slabs, i.e. $X_{12k}=X_{21k}$, $k=1,2$. The rank-1 approximation
problem is then
\be
\label{eq-prob3} {\rm min}_{{\bf y}\in\R^{2},{\bf z}\in\R^2}\;
||\ul{\bf X}-{\bf y}\otimes{\bf y}\otimes{\bf z}||^2\,. \ee

\noindent We can define a {\em symmetric slab rank} analogous to the
symmetric rank and propose a rank criterion similar to Lemma~\ref{lem-DSL}
and Lemma~\ref{lem-symrank}. Generic $2\times 2\times 2$ tensors with
symmetric slabs have {\emph{symmetric slab ranks} 2 and 3 on sets of positive
volume. Moreover, a result analogous to Theorem~\ref{t-1} and
Theorem~\ref{t-2} can be proven in this case.

\newpage
\section*{Appendix A: Proof of Theorem~\ref{t-1}}
\refstepcounter{section}
\setcounter{equation}{0}
\renewcommand{\thesection}{A}
\renewcommand{\theequation}{A.\arabic{equation}}
We make use of the first part of Section 2. Let $\ul{\bf X}$ be a generic
$2\times 2\times 2$ tensor with entries
\be
\label{eq-Xg}
\ul{\bf
X}=\left[\begin{array}{cc|cc}a&b&e&f\\c&d&g&h\end{array}\right]\,. \ee

\noindent We consider the rank-1 approximation problem (\ref{eq-prob}).
It is our goal to show that, for the optimal solution $\ul{\bf Y}={\bf
x}\otimes{\bf y}\otimes{\bf z}$, we have $\ul{\bf X}-\ul{\bf Y}$ in orbit
$D_3$. From the list of orbits in Table~\ref{tab-1}, it follows that it
suffices to show $\Delta( \ul{\bf X}-\ul{\bf Y})=0$ and $\ul{\bf X}-\ul{\bf
Y}$ has multilinear rank $(2,2,2)$. We will do this by considering the
stationary points of the rank-1 approximation problem. For later use, we
mention that the hyperdeterminant (\ref{eq-hypdet}) of $\ul{\bf X}$ in
(\ref{eq-Xg}) is given by
\be
\label{eq-hypdet2}
\Delta(\ul{\bf X})=(ah-bg+de-cf)^2-4\,(ad-bc)(eh-fg)\,.
\ee

\noindent We begin our proof by showing that for the best rank-1
approximation of $\ul{\bf X}$ we have $x_1\neq 0$, $x_2\neq 0$, $y_1\neq
0$, $y_2\neq 0$, $z_1\neq 0$ and $z_2\neq 0$ almost everywhere. Due to the
scaling indeterminacy in $({\bf x}\otimes{\bf y}\otimes{\bf z})$, this
implies that we may set $y_1=z_1=1$ without loss of generality.

\begin{lem}
\label{lem-nonzero}
Let $\ul{\bf X}$ be a generic $2\times 2\times 2$ tensor with a best
rank-$1$ approximation ${\bf x}\otimes{\bf y}\otimes{\bf z}$. Then $x_1\neq
0$, $x_2\neq 0$, $y_1\neq 0$, $y_2\neq 0$, $z_1\neq 0$, $z_2\neq 0$ almost
everywhere. \end{lem}

\noindent {\bf Proof.} We show that $z_1\neq 0$ almost everywhere. The
proofs for $y_1$, $y_2$ and $z_2$ are analogous. The proofs for $x_1$ and
$x_2$ follow by interchanging the roles of ${\bf x}$ and ${\bf y}$.

Let the criterion $\Psi$ be as in (\ref{eq-critsub}) and let $\Psi_0$
denote (\ref{eq-critsub}) with $z_1=0$. Then $\Psi<\Psi_0$ is equivalent to
\bdm
\hspace{-3cm}\left[(ay_1z_1+by_2z_1+ey_1z_2+fy_2z_2)^2+
(cy_1z_1+dy_2z_1+gy_1z_2+hy_2z_2)^2\right]\;>
\edm
\be
~\hspace{5cm}(z_1^2+z_2^2)\left[(ey_1+fy_2)^2+(gy_1+hy_2)^2\right]\,,
\ee

\noindent which, after setting $y_1=z_1=1$, can be rewritten as
\begin{eqnarray}
\label{eq-LL0}
(a^2+c^2-e^2-g^2+2z_2\,(ae+cg))+ &&\nonumber\\[2mm]
2y_2\,(ab+cd-ef-gh+z_2\,(af+be+ch+dg))+ && \nonumber\\[2mm]
y_2^2\,(b^2+d^2-f^2-h^2+2z_2\,(bf+dh)) &>&0\,.
\end{eqnarray}

\noindent Since $(bf+dh)\neq 0$ almost everywhere, it is possible to
choose $z_2$ such that the coefficient of $y_2^2$ is positive. Then there
is a range of values $y_2$ for which (\ref{eq-LL0}) holds. This shows that,
almost everywhere, we can find a better rank-1 approximation than setting
$z_1=0$. This completes the proof of $z_1\neq 0$.\ep

\noindent As mentioned above Lemma~\ref{lem-nonzero}, we set $y_1=z_1=1$
without loss of generality. Since the optimal ${\bf x}$ is
given by (\ref{eq-xopt}), the problem of finding a best rank-1
approximation of $\ul{\bf X}$ is now a problem in the variables $y_2$ and
$z_2$ only.

Next, we rewrite equations (\ref{eq-gradsub2y}) and (\ref{eq-gradsub2z})
specifying the stationary points $(y_2,z_2)$ as
\bdm
z_2^2\;\left[(ef+gh)\,y_2^2+(e^2+g^2-f^2-h^2)\,y_2-(ef+gh)\right]+
\edm
\bdm
z_2\;\left[(af+be+ch+dg)\,y_2^2+2\,(ae+cg-bf-dh)\,y_2
-(af+be+ch+dg)\right]+
\edm
\be
\label{eq-derLy22}
\left[(ab+cd)\,y_2^2+(a^2+c^2-b^2-d^2)\,y_2-(ab+cd)\right]=0\,,
\ee

\noindent and
\bdm
z_2^2\;\left[(bf+dh)\,y_2^2+(af+be+ch+dg)\,y_2+(ae+cg)\right]+
\edm
\bdm
z_2\;\left[(b^2+d^2-f^2-h^2)\,y_2^2+2\,(ab+cd-ef-gh)\,y_2
+(a^2+c^2-e^2-g^2)\right]+
\edm
\be
\label{eq-derLz22}
\left[-(bf+dh)\,y_2^2-(af+be+ch+dg)\,y_2-(ae+cg)\right]=0\,.
\ee

\noindent Using the expression (\ref{eq-xopt}) for ${\bf x}$, also the
hyperdeterminant (\ref{eq-hypdet}) of $\ul{\bf X}-\ul{\bf Y}$ can be
written as a function of $(y_2,z_2)$ only. After some manipulations, we
obtain
\begin{eqnarray}
\label{eq-discr2}
(1+y_2^2)^2(1+z_2^2)^2\;\;\Delta(\ul{\bf X}-\ul{\bf Y}) &=&
\left[z_2^2\;\left[(bg-de)\,y_2^2+(ag-bh-ce+df)\,y_2+(cf-ah)\right]\right.+
\nonumber \\[2mm]
&& z_2\;\left[(bc-ad+eh-fg)\,y_2^2+(bc-ad+eh-fg)\right]+ \nonumber \\[2mm]
&&\left.\left[(ah-cf)\,y_2^2+(ag-bh-ce+df)\,y_2+(de-bg)\right]\right]^2\,.
\end{eqnarray}

\noindent Equations (\ref{eq-derLy22}) and (\ref{eq-derLz22})
specifying the stationary points $(y_2,z_2)$, and the hyperdeterminant
(\ref{eq-discr2}) without the square, are of the same form: a
polynomial of degree 4 in $y_2$ and $z_2$ that is quadratic in both $y_2$
and $z_2$. We use the result of the following lemma to compare the
stationary points satisfying (\ref{eq-derLy22}) and (\ref{eq-derLz22}) to
the roots of (\ref{eq-discr2}).

\begin{lem}
\label{lem-root}
Let $f(u)=\alpha\,u^2+\beta\,u+\gamma$ and
$g(u)=\delta\,u^2+\epsilon\,u+\nu$ be second degree polynomials. Then $f$
and $g$ have a common root if and only if
\be
\label{eq-commonroot}
(\alpha\epsilon -\beta\delta)\,(\beta\nu -\epsilon\gamma)=
(\gamma\delta -\alpha\nu)^2\,.
\ee

\noindent Moreover, if $(\gamma\delta -\alpha\nu)$ and $(\alpha\epsilon
-\beta\delta)$ are nonzero, the common root is given by
\be
\label{eq-actualroot}
\frac{(\beta\nu -\epsilon\gamma)}{(\gamma\delta -\alpha\nu)}=
\frac{(\gamma\delta -\alpha\nu)}{(\alpha\epsilon -\beta\delta)}\,.
\ee
\end{lem}

\noindent {\bf Proof.} First, suppose $f$ and $g$ have a common root $r$.
Then $f(u)=\alpha(u-r)(u-r_1)$ and $g(u)=\delta(u-r)(u-r_2)$ for some $r_1$
and $r_2$. It follows that
\be
\beta=-\alpha\,(r+r_1)\quad\quad\epsilon=-\delta\,(r+r_2)\quad\quad
\gamma=\alpha\,r\,r_1\quad\quad\nu=\delta\,r\,r_2\,.
\ee

\noindent Using these expressions, it can be verified that
(\ref{eq-commonroot}) holds, and $r$ equals the expressions in
(\ref{eq-actualroot}).

Next, suppose (\ref{eq-commonroot}) holds. Let $f(u)=\alpha(u-r_1)(u-r_2)$
and $g(u)=\delta(u-r_3)(u-r_4)$ for some $r_1$, $r_2$, $r_3$, $r_4$. It
follows that
\be
\beta=-\alpha\,(r_1+r_2)\quad\quad\epsilon=-\delta\,(r_3+r_4)\quad\quad
\gamma=\alpha\,r_1\,r_2\quad\quad\nu=\delta\,r_3\,r_4\,.
\ee

\noindent Substituting these expressions into (\ref{eq-commonroot}) and
dividing both sides by $\alpha^2\delta^2$ yields
\be
(r_1+r_2-r_3-r_4)\,(r_1r_2\,(r_3+r_4)-r_3r_4\,(r_1+r_2))=
(r_1r_2-r_3r_4)^2\,.
\ee

\noindent This can be rewritten as
\be
(r_1-r_3)(r_1-r_4)(r_2-r_3)(r_2-r_4)=0\,,
\ee

\noindent which implies that $f$ and $g$ must have a common root. As above,
we have the expressions (\ref{eq-actualroot}) for the common root.
\ep

\noindent Using Lemma~\ref{lem-root}, the stationary points $(y_2,z_2)$
are found as follows. Equations (\ref{eq-derLy22})-(\ref{eq-derLz22})
represent two quadratic polynomials in $y_2$ that have a common root.
Lemma~\ref{lem-root} states that (\ref{eq-commonroot}) must hold, where all
coefficients are second degree polynomials in $z_2$. We rewrite this
equation as $P_z^{\rm stat}(z_2)=0$, where $P_z^{\rm stat}$ is a polynomial
of degree 8. The 8 roots of $P_z^{\rm stat}$ are the $z_2$ corresponding to
stationary points. For each $z_2$, the corresponding $y_2$ is the common
root given by (\ref{eq-actualroot}). Hence, there are 8 stationary points
$(y_2,z_2)$, and some of these may be complex.

Instead of interpreting (\ref{eq-derLy22})-(\ref{eq-derLz22}) as
polynomials in $y_2$, we may interpret them as polynomials in $z_2$ with
coefficients depending on $y_2$. As above, the $y_2$ of the stationary
points are then found by finding the roots of an 8th degree polynomial
$P_y^{\rm stat}(y_2)$ that is defined by (\ref{eq-commonroot}). For each
$y_2$, the corresponding $z_2$ is the common root given by
(\ref{eq-actualroot}). Both ways of obtaining the stationary points
necessarily yield the same result.

Analogously, we may determine the points $(y_2,z_2)$ satisfying
(\ref{eq-derLy22}) and having $\Delta(\ul{\bf X}-\ul{\bf Y})=0$ in
(\ref{eq-discr2}). The same approach yields the points satisfying
(\ref{eq-derLz22}) that are roots of (\ref{eq-discr2}). We denote
the 8th degree polynomials corresponding to (\ref{eq-derLy22}) and the
roots of (\ref{eq-discr2}) as $P_y^{\rm eig1}$ and $P_z^{\rm eig1}$.  We
denote the 8th degree polynomials corresponding to (\ref{eq-derLz22}) and
the roots of (\ref{eq-discr2}) as $P_y^{\rm eig2}$ and $P_z^{\rm eig2}$.
Using this approach, we obtain the following relation between the
stationary points and the roots of (\ref{eq-discr2}).

\newpage
\begin{table}[tp]
\begin{center}
\begin{tabular}{|c|cccccc|cc|}
\hline

(\ref{eq-derLy22}) and (\ref{eq-derLz22}) &

$y^{(1)}$ & $y^{(2)}$ & $y^{(3)}$ & $y^{(4)}$ & $y^{(5)}$ & $y^{(6)}$ &
$y^{(7)}$ & $y^{(8)}$ \\[2mm]

& $z^{(1)}$ & $z^{(2)}$ & $z^{(3)}$ & $z^{(4)}$ & $z^{(5)}$ & $z^{(6)}$ &
$z^{(7)}$ & $z^{(8)}$ \\

\hline

(\ref{eq-derLy22}) and root of (\ref{eq-discr2}) &

$y^{(1)}$ & $y^{(2)}$ & $y^{(3)}$ & $y^{(4)}$ & $y^{(5)}$ & $y^{(6)}$ &
$y^{(9)}$ & $y^{(10)}$ \\[2mm]

& $z^{(1)}$ & $z^{(2)}$ & $z^{(3)}$ & $z^{(4)}$ & $z^{(5)}$ & $z^{(6)}$ &
$z^{(7)}$ & $z^{(8)}$ \\

\hline

(\ref{eq-derLz22}) and root of (\ref{eq-discr2}) &

$y^{(1)}$ & $y^{(2)}$ & $y^{(3)}$ & $y^{(4)}$ & $y^{(5)}$ & $y^{(6)}$ &
$y^{(7)}$ & $y^{(8)}$ \\[2mm]

& $z^{(1)}$ & $z^{(2)}$ & $z^{(3)}$ & $z^{(4)}$ & $z^{(5)}$ & $z^{(6)}$ &
$z^{(9)}$ & $z^{(10)}$ \\

\hline

\end{tabular}
\end{center}

\caption{Schedule of points $(y_2,z_2)$ satisfying each pair of the
equations (\ref{eq-derLy22}), (\ref{eq-derLz22}), root of (\ref{eq-discr2}).
Equations (\ref{eq-derLy22})-(\ref{eq-derLz22}) describe stationary points,
while the roots of (\ref{eq-discr2}) have $\Delta(\ul{\bf X}-\ul{\bf Y})=0$.
As can be seen, the points $(y^{(i)},z^{(i)})$, $i=1,\ldots,6$, satisfy all
three equations.} \label{tab-points} \end{table}

\begin{lem}
\label{lem-points}
The points $(y_2,z_2)$ satisfying two of the three equations
$(\ref{eq-derLy22})$, $(\ref{eq-derLz22})$, root of $(\ref{eq-discr2})$,
are related as specified in Table~$\ref{tab-points}$. In particular, $6$ of
the $8$ stationary points are roots of $(\ref{eq-discr2})$. \end{lem}

\noindent {\bf Proof.} Using symbolic computation software, it can be
verified that
\be
\label{eq-bewijs1}
\frac{P_z^{\rm stat}(z_2)}{P_z^{\rm
eig1}(z_2)}=1\,,\quad\quad\quad\quad \frac{P_z^{\rm stat}(z_2)}{P_z^{\rm
eig2}(z_2)}= \frac{(eh-fg)\,z_2^2+(ah-bg+de-cf)\,z_2+(ad-bc)}
{(ad-bc)\,z_2^2-(ah-bg+de-cf)\,z_2+(eh-fg)}\,, \ee \be \label{eq-bewijs2}
\frac{P_y^{\rm stat}(y_2)}{P_y^{\rm eig2}(y_2)}=1\,,\quad\quad\quad\quad
\frac{P_y^{\rm stat}(y_2)}{P_y^{\rm eig1}(y_2)}=
\frac{(df-bh)\,y_2^2-(ah+bg-de-cf)\,y_2+(ce-ag)}
{(ce-ag)\,y_2^2+(ah+bg-de-cf)\,y_2+(df-bh)}\,. \ee

\noindent Hence, the roots $z_2$ of $P_z^{\rm stat}$ and $P_z^{\rm eig1}$
are identical, and so are the roots $y_2$ of $P_y^{\rm stat}$ and $P_y^{\rm
eig2}$. Also, $P_z^{\rm stat}$ and $P_z^{\rm eig2}$ have 6 of the 8 roots
in common, as do $P_y^{\rm stat}$ and $P_y^{\rm eig1}$. This implies that
the $z_2$-values of the stationary points coincide with the $z_2$-values of
the points satisfying (\ref{eq-derLy22}) that are roots of
(\ref{eq-discr2}). Analogously, the $y_2$-values of the stationary points
coincide with the $y_2$-values of the points satisfying (\ref{eq-derLz22})
that are roots of (\ref{eq-discr2}). Also, 6 of the $z_2$-values of
the stationary points coincide with the $z_2$-values of the points
satisfying (\ref{eq-derLz22}) that are roots of (\ref{eq-discr2}). And 6 of
the $y_2$-values of the stationary points coincide with the $y_2$-values of
the points satisfying (\ref{eq-derLy22}) that are roots of
(\ref{eq-discr2}).

In order to prove the relations in Table~\ref{tab-points}, it remains to
show that the 6 common $y_2$-values and the 6 common $z_2$-values form 6
common points $(y_2,z_2)$. Let $z_2$ be a root of $P_z^{\rm stat}$ and,
hence, of $P_z^{\rm eig1}$. The corresponding $y_2$ of the stationary point
is the common root given by (\ref{eq-actualroot}).
The corresponding $y_2$ of the point satisfying (\ref{eq-derLy22}) that is
a root of (\ref{eq-discr2}) is given by an analogous expression. Equating
these two expressions for $y_2$ yields an 8th degree polynomial in $z_2$
analogous to (\ref{eq-commonroot}). We denote this polynomial as $P_z^{\rm
com}$. Using symbolic computation software, it can be verified that
\be
\label{eq-bewijs3}
\frac{P_z^{\rm stat}(z_2)}{P_z^{\rm
com}(z_2)}= \frac{(eh-fg)\,z_2^2+(ah-bg+de-cf)\,z_2+(ad-bc)}
{(ef+gh)\,z_2^2+(af+be+ch+dg)\,z_2+(ab+cd)}\,. \ee

\noindent Hence, $P_z^{\rm stat}$ and $P_z^{\rm com}$ have 6 common roots.
This implies that 6 stationary points $(y_2,z_2)$ are also roots of
(\ref{eq-discr2}). This completes the proof of the relations in
Table~\ref{tab-points}.\ep

\noindent So far, we have shown that 6 of the 8 stationary points in the
rank-1 approximation problem satisfy $\Delta(\ul{\bf X}-\ul{\bf Y})=0$. In
Lemma~\ref{lem-2points} below, we show that the two other stationary
points $(y^{(7)},z^{(7)})$ and $(y^{(8)},z^{(8)})$ correspond to ${\bf
x}={\bf 0}$ in (\ref{eq-xopt}), which is not a best rank-1 approximation.
The global minimum of the rank-1 approximation problem is thus attained in
one of the stationary points $(y^{(i)},z^{(i)})$, $i=1,\ldots,6$. In
Lemma~\ref{lem-m222} the proof of Theorem~\ref{t-1} is completed by showing
that the multilinear rank of $\ul{\bf X}-\ul{\bf Y}$ equals $(2,2,2)$ for
these stationary points. Together with $\Delta(\ul{\bf X}-\ul{\bf Y})=0$,
this implies that $\ul{\bf X}-\ul{\bf Y}$ is in orbit $D_3$.

Next, we consider the two stationary points $(y^{(7)},z^{(7)})$
and $(y^{(8)},z^{(8)})$. Note that
$y^{(7)}$ and $y^{(8)}$ are the roots of the numerator of
(\ref{eq-bewijs2}), $y^{(9)}$ and $y^{(10)}$ are the roots of the
denominator of (\ref{eq-bewijs2}), $z^{(7)}$ and $z^{(8)}$ are the roots of
the numerator of (\ref{eq-bewijs1}), and $z^{(9)}$ and $z^{(10)}$ are the
roots of the denominator of (\ref{eq-bewijs1}). Moreover, these four
polynomials of degree 2 have identical discriminant that is equal to the
hyperdeterminant of $\ul{\bf X}$ as given in (\ref{eq-hypdet2}).

Hence, if $\Delta(\ul{\bf X})<0$, i.e. $\ul{\bf X}$ is in orbit $G_3$,
then the stationary points $(y^{(7)},z^{(7)})$ and $(y^{(8)},z^{(8)})$
are complex. Since we only consider real-valued rank-1 approximations, we
discard these two stationary points. If $\Delta(\ul{\bf X})>0$, i.e.
$\ul{\bf X}$ is in orbit $G_2$, we resort to Lemma~\ref{lem-2points}.

\begin{lem}
\label{lem-2points}
Suppose $\Delta(\ul{\bf X})>0$. Then the stationary points
$(y^{(7)},z^{(7)})$ and $(y^{(8)},z^{(8)})$ in Table~$\ref{tab-points}$
yield ${\bf x}={\bf 0}$ in $(\ref{eq-xopt})$, and do not correspond to the
global minimum almost everywhere.\end{lem}

\noindent {\bf Proof.} It can be verified that $y^{(7)}$ and $y^{(8)}$ are
given by
\be
\label{eq-staty}
\frac{(ah+bg-de-cf)\pm\sqrt{(ah+bg-de-cf)^2-4(df-bh)(ce-ag)}}{2\,
(df-bh)}\,,
\ee

\noindent and $z^{(7)}$ and $z^{(8)}$ are given by
\be
\label{eq-statz}
\frac{-(ah-bg+de-cf)\pm\sqrt{(ah-bg+de-cf)^2-4(eh-fg)(ad-bc)}}{2\,
(eh-fg)}\,,
\ee

\noindent where $\pm$ is $+$ in one stationary point and $-$ in the other.
Using symbolic computation software, it can be verified that the expression
for ${\bf x}$ in (\ref{eq-xopt}) is all-zero for $(y^{(7)},z^{(7)})$ and
$(y^{(8)},z^{(8)})$. Hence, both stationary points yield the all-zero
solution. This is not the global minimum since the solution
\be
{\bf x}=\left(\begin{array}{c} a\\
0\end{array}\right)\,,\quad\quad\quad
 {\bf y}=\left(\begin{array}{c}1\\
0\end{array}\right)\,,\quad\quad\quad
{\bf z}=\left(\begin{array}{c}1\\
0\end{array}\right)\,,
\ee

\noindent yields a lower criterion $\Psi$ in (\ref{eq-crit}) when
$a\neq 0$. This completes the proof.\ep

\begin{lem}
\label{lem-m222}
For the stationary points $(y^{(i)},z^{(i)})$, $i=1,\ldots,6$, in
Table~$\ref{tab-points}$ the multilinear rank of $\ul{\bf X}-\ul{\bf Y}$
equals $(2,2,2)$ almost everywhere. \end{lem}

\noindent {\bf Proof.} Let $\ul{\bf Z}=\ul{\bf X}-\ul{\bf Y}=\ul{\bf
X}-{\bf x}\otimes{\bf y}\otimes{\bf z}$, where ${\bf x}$ is given by
(\ref{eq-xopt}), $y_1=z_1=1$, and $(y_2,z_2)$ is a stationary point. If one
of the frontal slabs ${\bf Z}_1$ and ${\bf Z}_2$ of $\ul{\bf Z}$ is
nonsingular, then the mode-1 and mode-2 ranks of $\ul{\bf Z}$ are equal to
2. Next, we show that det$({\bf Z}_1)=$ det$({\bf Z}_2)=0$ corresponds to a
set of measure zero. It can be verified that \be {\rm det}({\bf
Z}_1)=\frac{z_2\,[-(de-bg)-(ag-bh-ce+df)\,y_2-
(ah-cf)\,y_2^2+(ad-bc)(1+y_2^2)\,z_2]}{(1+y_2^2)(1+z_2^2)}\,, \ee

\noindent and
\be
{\rm det}({\bf Z}_2)=\frac{(eh-fg)(1+y_2^2)+z_2\,
[-(ah-cf)+(ag-bh-ce+df)\,y_2-(de-bg)\,y_2^2]}{(1+y_2^2)(1+z_2^2)}\,.
\ee

\noindent Suppose det$({\bf Z}_1)=$ det$({\bf Z}_2)=0$, i.e. the numerators
of the above expressions are zero. Since $z_2\neq 0$ almost everywhere (see
Lemma~\ref{lem-nonzero}), we divide the numerator of det$({\bf Z}_1)$
by $z_2$. We then obtain two equations of the form $z_2=s(y_2)/t(y_2)$.
Equating both expressions for $z_2$ yields a fourth degree polynomial in
$y_2$ that can be written as
\be
[(ag-ce)\,y_2^2-(ah+bg-cf-de)\,y_2+(bh-df)]\;
[(df-bh)\,y_2^2-(ah+bg-cf-de)\,y_2+(ce-ag)]=0\,.
\ee

\noindent These two second degree polynomials are the numerator (times
$-1$) and denominator of (\ref{eq-bewijs2}). As explained above, the roots
of these polynomials are complex if $\Delta(\ul{\bf X})<0$. In this case,
it is not possible to choose $y_2$ and $z_2$ such that det$({\bf Z}_1)=$
det$({\bf Z}_2)=0$. When $\Delta(\ul{\bf X})>0$, the sought values of $y_2$ are
$y^{(7)}$, $y^{(8)}$, $y^{(9)}$ and $y^{(10)}$. Therefore, in this case we
may conclude that the points $(y_2,z_2)$ for which det$({\bf Z}_1)=$
det$({\bf Z}_2)=0$ are not among the first 6 stationary points in
Table~\ref{tab-points} almost everywhere.

Hence, the multilinear rank of $\ul{\bf Z}$ equals $(2,2,*)$. If one of the
top and bottom slabs of $\ul{\bf Z}$ is nonsingular, then also its mode-3
rank equals 2. A proof of this can be obtained analogous as above by
interchanging the roles of ${\bf x}$ and ${\bf z}$. This completes the
proof. \ep

\newpage
\subsection*{Numerical examples}
Here, we illustrate the proof of Theorem~\ref{t-1} by means of two
examples. We take two random $\ul{\bf X}$, one that has
$\Delta(\ul{\bf X})>0$ (orbit $G_2$) and one that has $\Delta(\ul{\bf
X})<0$ (orbit $G_3$).

Our first example is
\be
\ul{\bf X}=\left[\begin{array}{cc|cc}
-0.4326  &  0.1253 & -1.1465 &   1.1892 \\
-1.6656  &  0.2877 &  1.1909 &  -0.0376
\end{array}\right]\,. \ee

\noindent We have $\Delta(\ul{\bf X})=2.7668$. In
the table below, we list the stationary points $(y_2,z_2)$, their values of
$\Psi$ in (\ref{eq-critsub}), their values of $\Delta(\ul{\bf X}-\ul{\bf
Y})$, and state whether their Hessian matrix is positive definite or not.
Two of the stationary points $(y^{(i)},z^{(i)})$, $i=1,\ldots,6$, are
complex. The remaining four points are the first four points in the table,
and have $\Delta(\ul{\bf X}-\ul{\bf Y})$ close to zero. The second point
corresponds to the global minimum and is also found when computing a best
rank-1 approximation to $\ul{\bf X}$ via an alternating least squares
algorithm. For $\ul{\bf Z}=\ul{\bf X}-\ul{\bf Y}$, the matrix ${\bf
Z}_2{\bf Z}_1^{-1}$ has a double eigenvalue $0.9185$ with only one
associated eigenvector. Lemma~\ref{lem-DSL} implies that $\ul{\bf
Z}$ is in orbit $D_3$. The last two points in the table are the stationary
points $(y^{(7)},z^{(7)})$ and $(y^{(8)},z^{(8)})$. From
Lemma~\ref{lem-2points} it follows that they have $\Delta(\ul{\bf
X}-\ul{\bf Y})=\Delta(\ul{\bf X})$ and $\Psi=||\ul{\bf X}||^2$.
\vspace{5mm}

\begin{center}
\begin{tabular}{|cc|ccc|}
\hline

$y_2$ & $z_2$ & $\Psi$ & $\Delta(\ul{\bf X}-\ul{\bf Y})$ & Hessian PD
\\[2mm]
\hline
-0.592958 & 0.621735 &  5.1164 & 1.4166e-12 & no \\[2mm]
-0.229249 &-1.08855  &  2.6863 & 9.6802e-13 & yes \\[2mm]
2.22613  &  0.452035 &  7.1313 & 2.1210e-12 & no \\[2mm]
2.42488  & -2.88759 &  6.5289 & 1.2999e-14 & no\\[2mm]
\hline
1.17156  & 1.15843  &  7.2081 & 2.7668  & no \\[2mm]
5.96728  & -0.05296 &  7.2081 & 2.7668  & no\\

\hline

\end{tabular}
\end{center}
\vspace{5mm}

\noindent Our second example is
\be
\ul{\bf X}=\left[\begin{array}{cc|cc}
-1.6041 &  -1.0565 & 0.8156 &   1.2902 \\
 0.2573 &   1.4151 & 0.7119 &   0.6686
\end{array}\right]\,. \ee

\noindent We have $\Delta(\ul{\bf X})=-2.7309$.
In the table below, we list the stationary points $(y_2,z_2)$ in the same
way as in the first example. Two of the stationary points
$(y^{(i)},z^{(i)})$, $i=1,\ldots,6$, are complex. Since $\Delta(\ul{\bf
X})<0$, the points $(y^{(7)},z^{(7)})$ and $(y^{(8)},z^{(8)})$ are
also complex. Hence, four real stationary points are left, that all have
$\Delta(\ul{\bf X}-\ul{\bf Y})$ close to zero. The first point in the table
corresponds to the global minimum and is also found when computing a best
rank-1 approximation to $\ul{\bf X}$ via an alternating least squares
algorithm. For $\ul{\bf Z}=\ul{\bf X}-\ul{\bf Y}$, the matrix
${\bf Z}_2{\bf Z}_1^{-1}$ has a double eigenvalue
$1.6712$ with only one associated eigenvector. Lemma~\ref{lem-DSL} implies
that $\ul{\bf Z}$ is in orbit $D_3$.\vspace{5mm}

\begin{center}
\begin{tabular}{|cc|ccc|}
\hline

$y_2$ & $z_2$ & $\Psi$ & $\Delta(\ul{\bf X}-\ul{\bf Y})$ & Hessian PD \\[2mm]
\hline
0.995675 & -0.598339 & 3.1185 & 1.3801e-11 & yes \\[2mm]
-0.865475 & 0.0601889 & 8.2319  & 1.5479e-13 & no \\[2mm]
 2.06437 & 1.78102 & 6.6050  & 1.6050e-13 & no \\[2mm]
-0.675154 &  9.24487 & 9.0028 & 2.6216e-13 & no \\
\hline
\end{tabular}
\end{center}
\vspace{5mm}

\section*{Appendix B: Proof of Theorem~\ref{t-2}}
\refstepcounter{section}
\setcounter{equation}{0}
\renewcommand{\thesection}{B}
\renewcommand{\theequation}{B.\arabic{equation}}
We make use of the derivations in Section 5. Let $\ul{\bf X}$ be a generic
symmetric $2\times 2\times 2$ tensor (\ref{eq-Xgsym}). We consider the
symmetric rank-1 approximation problem (\ref{eq-prob2}). It is our goal to
show that, for the optimal solution $\ul{\bf Y}={\bf y}\otimes{\bf
y}\otimes{\bf y}$, we have $\ul{\bf X}-\ul{\bf Y}$ in orbit $D_3$. From the
list of orbits in Table~\ref{tab-orbit2}, it follows that it suffices to
show $\Delta( \ul{\bf X}-\ul{\bf Y})=0$ and $\ul{\bf X}-\ul{\bf Y}$ has
multilinear rank $(2,2,2)$. We will do this by considering the stationary
points of the symmetric rank-1 approximation problem.

Let $\ul{\bf Y}$ be as in (\ref{eq-Ysym}). The stationary points are given
by (\ref{eq-grad2}), which can be written as
\begin{eqnarray}
\label{eq-staty1}
y_1^5+y_1y_2^4+2\,y_1^3y_2^2-2\,b\,y_1y_2-a\,y_1^2-c\,y_2^2
&=& 0 \,, \\[2mm]
\label{eq-staty2}
y_2^5+y_1^4y_2+2\,y_1^2y_2^3-2\,c\,y_1y_2-d\,y_2^2-b\,y_1^2
&=& 0\,.
\end{eqnarray}

\noindent Note that the entries $a,b,c,d$ are nonzero almost everywhere.
If one of $y_1$ and $y_2$ is zero, it follows that both are
zero almost everywhere. Since this corresponds to an all-zero $\ul{\bf Y}$,
which is not the optimal solution, we may assume that $y_1\neq 0$ and
$y_2\neq 0$ almost everywhere.

Multiplying (\ref{eq-staty2}) by $y_1$ and subtracting $y_2$ times
(\ref{eq-staty1}) yields
\be
\label{eq-poly1}
-b\,y_1^3+(a-2c)\,y_1^2y_2+(2b-d)\,y_1y_2^2+c\,y_2^3=0\,.
\ee

\noindent Defining $z=y_1/y_2$ and dividing (\ref{eq-poly1}) by $y_2^3$,
we obtain
\be
\label{eq-poly2}
-b\,z^3+(a-2c)\,z^2+(2b-d)\,z+c=0\,.
\ee

\noindent This yields three solutions for $z=y_1/y_2$, two of which may be
complex. For each solution $z$, the corresponding stationary point
$(y_1,y_2)$ satisfying (\ref{eq-staty1})-(\ref{eq-staty2}) is given by
\be
\label{eq-staty1y2}
y_1=z\,y_2\,,\quad\quad\quad\quad
y_2^3=\frac{a\,z^2+2\,b\,z+c}{z^5+2\,z^3+z}=
\frac{b\,z^2+2\,c\,z+d}{z^4+2\,z^2+1}\,,
\ee

\noindent where the latter equality is equivalent to (\ref{eq-poly2}). The
polynomial (\ref{eq-poly2}) determining the stationary points is also
reported by \cite[section 3.5]{LMVrank1}.

Next, we consider the hyperdeterminant $\Delta(\ul{\bf X}-\ul{\bf Y})$.
For $y_1=z\,y_2$, we have
\be
\label{eq-XminY}
\ul{\bf X}-\ul{\bf Y}=\left[\begin{array}{cc|cc}
a-z^3\,y_2^3 & b-z^2\,y_2^3 & b-z^2\,y_2^3 & c-z\,y_2^3 \\
b-z^2\,y_2^3 & c-z\,y_2^3 & c-z\,y_2^3 & d-y_2^3\end{array}\right]\,.
\ee

\noindent Using (\ref{eq-hypdetsym}), we obtain
\be
\label{eq-hypdetpoly}
\Delta(\ul{\bf X}-\ul{\bf Y})=\Delta(\ul{\bf X})+f(z)\,y_2^3
+(a-3\,b\,z+3\,c\,z^2-d\,z^3)^2\,y_2^6\,,
\ee

\noindent with
\be
f(z)=[-4\,b^3+6\,abc-2\,a^2d]+z\,[6\,b^2c-12\,ac^2+6\,abd]
+z^2\,[6\,bc^2-12\,b^2d+6\,acd]+z^3\,[6\,bcd-2\,ad^2]\,.
\ee

\noindent We substitute the second expression for $y_2^3$ in
(\ref{eq-staty1y2}) into (\ref{eq-hypdetpoly}) and multiply by
$(z^4+2\,z^2+1)^2$. Using symbolic computation software, it can be verified
that this yields
\be
\label{eq-iszero}
(-b\,z^3+(a-2c)\,z^2+(2b-d)\,z+c)\,P(z)\,,
\ee

\noindent where $P(z)$ is a $7$th degree polynomial in $z$. By
(\ref{eq-poly2}), the expression (\ref{eq-iszero}) is identical to zero.
Hence, for all three stationary points $(y_1,y_2)$, we have $\Delta(\ul{\bf
X}-\ul{\bf Y})=0$.

In the final part of the proof, we show that $\ul{\bf X}-\ul{\bf Y}$ has
multilinear rank $(2,2,2)$ almost everywhere. Since the mode-$n$ ranks of
symmetric tensors are equal for each mode, it suffices to show that the two
slabs of (\ref{eq-XminY}) are nonsingular almost everywhere. Let $\ul{\bf
Z}=\ul{\bf X}-\ul{\bf Y}$. We have
\begin{eqnarray}
{\rm det}({\bf Z}_1) &=& (ac-b^2)+y_2^3\,(-c\,z^3+2\,b\,z^2-a\,z)\,, \\[2mm]
{\rm det}({\bf Z}_2) &=& (bd-c^2)+y_2^3\,(-d\,z^2+2\,c\,z-b)\,.
\end{eqnarray}

\noindent Hence, det$({\bf Z}_1)=$ det$({\bf Z}_2)=0$ implies
\be
(bd-c^2)(c\,z^3-2\,b\,z^2+a\,z)+(ac-b^2)(d\,z^2-2\,c\,z+b)=0\,,
\ee

\noindent which can be written as
\be
\label{eq-poly3}
z^3\,[bcd-c^3]+z^2\,[2\,bc^2+acd-3\,b^2d]+z\,[2\,b^2c+abd-3\,ac^2]+[abc-b^3]
= 0\,. \ee

\noindent Since the $3$rd degree polynomials (\ref{eq-poly2}) and
(\ref{eq-poly3}) do not have generically common roots, it
follows that at least one of the slabs ${\bf Z}_1$ and ${\bf Z}_2$ is
nonsingular almost everywhere. As explained above, this implies that
$\ul{\bf X}-\ul{\bf Y}$ has multilinear rank $(2,2,2)$ almost everywhere.
This completes the proof of Theorem~\ref{t-2}.

\section*{Appendix C: Orbits $D_3$ and $G_3$ of real symmetric $2\times
2\times 2$ tensors}
\refstepcounter{section}
\setcounter{equation}{0}
\renewcommand{\thesection}{C}
\renewcommand{\theequation}{C.\arabic{equation}}
Here, we show that any real symmetric $2\times 2\times 2$ tensor $\ul{\bf
X}$ in orbit $D_3$ or $G_3$ is related to the canonical form $\ul{\bf Y}$
of the orbit by an invertible multilinear transformation $({\bf S},{\bf
S},{\bf S})\cdot\ul{\bf Y}=\ul{\bf X}$.

First, we consider orbit $D_3$, which is defined by symmetric rank 3,
multilinear rank $(2,2,2)$, and hyperdeterminant $\Delta=0$. It follows
from the proof of Proposition~\ref{p-orbit2} that we may assume without
loss of generality that $\ul{\bf X}$ in orbit $D_3$ has the form
\be
\label{eq-oftheform2}
\ul{\bf X}=\left[\begin{array}{cc|cc} a & 1 & 1 & 1\\
1 & 1 & 1 & d\end{array}\right]\,, \ee

\noindent with
\be
\label{eq-hypdetzero}
\Delta(\ul{\bf X})=a^2d^2-6ad+4a+4d-3=0\,.
\ee

\noindent Our goal is to find a nonsingular
\be
{\bf S}=\left[\begin{array}{cc}
s_1 & s_2 \\ s_3 & s_4\end{array}\right]\,,
\ee

\noindent such that $({\bf S},{\bf
S},{\bf S})\cdot\ul{\bf Y}=\ul{\bf X}$, where the canonical form $\ul{\bf
Y}$ of orbit $D_3$ is given in Table~\ref{tab-orbit2}, i.e. \be
\label{eq-SX} ({\bf S},{\bf S},{\bf S})\cdot\;\left[\begin{array}{cc|cc} 0
& 1 & 1 & 0\\ 1 & 0 & 0 &
0\end{array}\right]=\left[\begin{array}{cc|cc}a&1&1&1\\ 1 & 1 & 1 &
d\end{array}\right]\,. \ee

\noindent This yields the following four equations:
\be
\label{eq-D31}
3\,s_1^2s_2=a\,,\quad\quad\quad\quad 3\,s_3^2s_4=d\,,
\ee
\be
\label{eq-D32}
s_1^2s_4+2\,s_1s_2s_3=1\,,\quad\quad\quad\quad
s_2s_3^2+2\,s_1s_3s_4=1\,.
\ee

\noindent Note that the case $a=d=1$ has $\Delta=0$ but yields multilinear
rank $(1,1,1)$ and, hence, is not included in orbit $D_3$. The case $a=0$,
$d=3/4$ is in orbit $D_3$ and its solution of (\ref{eq-D31})-(\ref{eq-D32})
is
\be
\label{eq-S}
{\bf S}=\left[\begin{array}{cc}
1&0\\ 1/2 & 1\end{array}\right]\,.
\ee

\noindent The case $a=3/4$, $d=0$ can be treated analogously. In the
remaining part of the proof we assume $a\neq 0$ and $d\neq 0$. This implies
that all entries of ${\bf S}$ are nonzero. From (\ref{eq-hypdetzero}) it
follows that
\be
\label{eq-hypdetzerod}
d=\frac{3\,a-2\pm 2\,(1-a)\sqrt{1-a}}{a^2}\,.
\ee

\noindent Hence, we must have $a<1$. Since (\ref{eq-hypdetzero}) is
symmetric in $a$ and $d$, also $d<1$ must hold.

Next, we solve the system (\ref{eq-D31})-(\ref{eq-D32}).
From (\ref{eq-D31}) we get $s_2=a/(3\,s_1^2)$ and $s_4=d/(3\,s_3^2)$.
Substituting this into (\ref{eq-D32}) yields, after rewriting,
\be
\label{eq-D33}
\frac{d}{3}\;\left(\frac{s_1}{s_3}\right)^3=\left(\frac{s_1}{s_3}\right)\;-
\;\frac{2\,a}{3}\,,\quad\quad\quad\quad
\frac{d}{3}\;\left(\frac{s_1}{s_3}\right)^3=\frac{1}{2}\;\left(\frac{s_1}{s_
3}\right)^2\;-\;\frac{a}{6}\,.
\ee

\noindent We equate the right-hand sides of (\ref{eq-D33}), which yields
\be
\label{eq-sol13}
\left(\frac{s_1}{s_3}\right)=1\pm\sqrt{1-a}\,.
\ee

\noindent Substituting this into one equation of (\ref{eq-D33}) gives us
\be
d=\frac{3-2\,a\pm 3\sqrt{1-a}}{(1\pm\sqrt{1-a})^3}\,.
\ee

\noindent It can be verified that this expression for $d$ is identical to
(\ref{eq-hypdetzerod}). Hence, equation (\ref{eq-sol13}), together with
$s_2=a/(3\,s_1^2)$ and $s_4=d/(3\,s_3^2)$, solves the system
(\ref{eq-D31})-(\ref{eq-D32}). Note that since both tensors in
(\ref{eq-SX}) have multilinear rank $(2,2,2)$, it follows that ${\bf S}$
is nonsingular. Hence, we have shown that for any $\ul{\bf X}$ in orbit
$D_3$ there exists a nonsingular ${\bf S}$ such that $({\bf S},{\bf S},{\bf
S})\cdot\ul{\bf Y}=\ul{\bf X}$, where $\ul{\bf Y}$ is the canonical form of
orbit $D_3$.

Next, we consider orbit $G_3$, which is defined by symmetric rank 3,
multilinear rank $(2,2,2)$, and hyperdeterminant $\Delta<0$. As above, we
may assume that $\ul{\bf X}$ in $G_3$ has the form (\ref{eq-oftheform2})
with
\be
\label{eq-hypdetneg}
\Delta(\ul{\bf X})=a^2d^2-6ad+4a+4d-3<0\,.
\ee

\noindent It is our goal to find nonsingular ${\bf S}$ in (\ref{eq-S}) such
that
\be
\label{eq-SX2}
({\bf S},{\bf S},{\bf S})\cdot\;\left[\begin{array}{cc|cc} -1 & 0 & 0 & 1\\
0 & 1 & 1 & 0\end{array}\right]=\left[\begin{array}{cc|cc}a&1&1&1\\
1 & 1 & 1 & d\end{array}\right]\,,
\ee

\noindent where the former tensor is the canonical form of orbit
$G_3$ as given in Table~\ref{tab-orbit2}. This yields the following four
equations: \be \label{eq-G31} -s_1^3+3\,s_1s_2^2=a\,,\quad\quad\quad\quad
-s_3^3+3\,s_3s_4^2=d\,, \ee \be \label{eq-G32}
-s_1^2s_3+2\,s_1s_2s_4+s_2^2s_3=1\,,\quad\quad\quad\quad
-s_1s_3^2+2\,s_2s_3s_4+s_1s_4^2=1\,. \ee

\noindent The case $a=0$, $d<3/4$ has solution $s_1=0$, $s_3^3=3/4-d>0$,
$s_2^2=1/s_3$, $s_4^2=1/(4\,s_3)$, with det$({\bf S})=-\sqrt{s_3}<0$. The
case $a<3/4$, $d=0$ can be treated analogously. In the remaining part of
the proof we assume $a\neq 0$ and $d\neq 0$. This implies that $s_1$ and
$s_3$ are nonzero. Note that $a=1$ implies
$\Delta=(d-1)^2$, which is not in orbit $G_3$. Analogously, $d=1$
is not in orbit $G_3$ either. In fact, $\Delta<0$ implies $a<1$ and $d<1$.

Next, we solve the system (\ref{eq-G31})-(\ref{eq-G32}).
Expressions for $s_2$ and $s_4$ are obtained from (\ref{eq-G31}) as
\be
\label{eq-sol24}
s_2^2=\frac{s_1^3+a}{3\,s_1}\,,\quad\quad\quad\quad
s_4^2=\frac{s_3^3+d}{3\,s_3}\,.
\ee

\noindent Equations (\ref{eq-G32}) can be written as
\be
2\,s_2s_4=\frac{1+s_1^2s_3-s_2^2s_3}{s_1}\,,\quad\quad\quad\quad
2\,s_2s_4=\frac{1+s_1s_3^2-s_1s_4^2}{s_3}\,.
\ee

\noindent Equating the right-hand sides and substituting (\ref{eq-sol24})
yields, after rewriting,
\be
\label{eq-poly13}
d\;\left(\frac{s_1}{s_3}\right)^3-3\;\left(\frac{s_1}{s_3}\right)^2+3\;\left
(\frac{s_1}{s_3}\right)-a=0\,.
\ee

\noindent The discriminant of this 3rd degree polynomial equals
$-27\Delta>0$, which implies that (\ref{eq-poly13}) has three distinct real
roots. Let $s_1=\alpha\,s_3$, where the root $\alpha$ satisfies
\be
\label{eq-alpha}
d\,\alpha^3=3\,\alpha^2-3\,\alpha+a\,.
\ee

\noindent Substituting $s_1=\alpha\,s_3$ and (\ref{eq-sol24}) into the first
equation of (\ref{eq-G32}) yields
\be
4\,s_3^3\,(d\,\alpha^4+2a\,\alpha-3\,\alpha^2)=9-6a/\alpha+a^2/\alpha^2-
4ad \,\alpha\,. \ee

\noindent Using (\ref{eq-alpha}), this can be rewritten as
\be
\label{eq-sol3}
s_3^3=\frac{\alpha\,(3-d\,\alpha)^2-4ad}{12\,(\alpha^2-2\,\alpha+a)}\,.
\ee

\noindent It remains to verify that the expressions (\ref{eq-sol24}) are
nonnegative. Our proof is tedious and long. Below, we give a
summary of it. The full proof is available on request.

Substituting (\ref{eq-sol3}) and using (\ref{eq-alpha}), it can
be shown that the expressions (\ref{eq-sol24}) are nonnegative if
\be
\label{eq-P12}
P_1(\alpha)=\alpha^2\,(4-3d)+\alpha\,(ad-3)+a\ge 0\,,
\quad\quad\quad\quad
P_2(\alpha)=-d\,\alpha^2+\alpha\,(3-ad)-a\ge 0\,.
\ee

\noindent Note that $P_1+P_2=4\,\alpha^2\,(1-d)>0$. Also, the
leading coefficient of $P_1$ is always positive. The roots of $P_1$ are
given by
\be
r_1=\frac{3-ad-\sqrt{a^2d^2+6ad-16a+9}}{2\,(4-3d)}\,,\quad\quad\quad\quad
r_2=\frac{3-ad+\sqrt{a^2d^2+6ad-16a+9}}{2\,(4-3d)}\,.
\ee

\noindent The roots of $P_2$ are given by
\be
r_3=\frac{ad-3-\sqrt{a^2d^2-10ad+9}}{-2d}\,,\quad\quad\quad\quad
r_4=\frac{ad-3+\sqrt{a^2d^2-10ad+9}}{-2d}\,.
\ee

\noindent Let $P_3(x)=d\,x^3-3\,x^2+3\,x-a$. To prove (\ref{eq-P12}), we
focus on the sign of $P_3$ in the roots $r_1,r_2,r_3,r_4$. When the
discriminant of $P_1$ is nonnegative, we can distinguish three cases. In
these cases, the sign of $P_3$ in the roots $r_1$ and $r_2$ is as follows:
\begin{eqnarray}
\label{eq-case1}
{\rm case\;I}:\quad  a>0\quad{\rm and}\quad a^2d>-3a+4a\sqrt{a}
&\Rightarrow & P_3(r_1)\ge 0\quad P_3(r_2)\ge 0\,,\\[2mm] \label{eq-case2}
{\rm case\;II}:\quad  a>0\quad{\rm and}\quad a^2d<-3a-4a\sqrt{a} &
\Rightarrow & P_3(r_1)\le 0\quad P_3(r_2)\le 0\,, \\[2mm] \label{eq-case3}
{\rm case\;III}:\quad   a<0 & \Rightarrow & P_3(r_1)\ge 0\quad P_3(r_2)\ge
0\,. \end{eqnarray}

\noindent When the discriminant of $P_2$ is nonnegative, we have
\be
\label{eq-rootsP2}
P_3(r_4)\le 0\,,\quad\quad
P_3(r_3)\left\{\begin{array}{ll}
\ge 0 & {\rm if\;} d>0\,, \\
\le 0 & {\rm if\;} d<0\,.
\end{array}\right.
\ee

\noindent Suppose $d>0$. Then the leading coefficient of $P_2$ is negative
and its discriminant is positive (since $a<1$ and $d<1$). Hence, $P_2$ has
real roots. Recall that the leading coefficient of $P_1$ is always positive.
Suppose the roots of $P_1$ are real. Then we are in case I or
case III (since case II implies $d<0$). Since $P_1+P_2>0$, there must hold
$r_4\le r_1\le r_2\le r_3$. From (\ref{eq-case1}), (\ref{eq-case3}), and
(\ref{eq-rootsP2}), it follows that $P_3$ has a root $\alpha$ in the
interval $[r_4,r_1]$ for which (\ref{eq-P12}) holds. If the roots of $P_1$
are not real, then (\ref{eq-rootsP2}) implies that $P_3$ has a root
$\alpha$ in the interval $[r_4,r_3]$ for which (\ref{eq-P12}) holds.

Suppose next that $d<0$. Then the leading coefficients of $P_1$ and $P_2$
are positive. Suppose $P_1$ and $P_2$ both have real roots. Since
$P_1+P_2>0$, there must hold either $r_3\le r_4\le r_1\le r_2$ or $r_1\le
r_2\le r_3\le r_4$. From (\ref{eq-case1})-(\ref{eq-rootsP2}) we obtain the
following. Suppose we are in case I or case III. If $r_3\le r_4\le r_1\le
r_2$, then $P_3$ has a root $\alpha$ in the interval $[r_4,r_1]$ for which
(\ref{eq-P12}) holds. If $r_1\le r_2\le r_3\le r_4$, then $P_3$ has a root
$\alpha$ in the interval $[r_2,r_3]$ for which (\ref{eq-P12}) holds.
Suppose we are in case II. Then $P_3(r_3)\le 0$ and $P_3(r_1)\le 0$. From
the shape of $P_3$ it follows that it has a root $\alpha\le r_3$ if $r_3\le
r_4\le r_1\le r_2$, or a root $\alpha\le r_1$ if $r_1\le r_2\le r_3\le
r_4$. In both situations, we have (\ref{eq-P12}) for this root $\alpha$.

When $d<0$ and $P_1$ does not have real roots, (\ref{eq-rootsP2}) implies
that $P_3$ has a root $\alpha\le r_3$ for which (\ref{eq-P12}) holds. When
$d<0$ and $P_2$ does not have real roots, (\ref{eq-case1})-(\ref{eq-case3})
imply that $P_3$ cannot have all three roots in the interval $[r_1,r_2]$.
Hence, there exists a root $\alpha$ for which (\ref{eq-P12}) holds.
Finally, it can be shown that $P_1$ and $P_2$ cannot both have complex
roots when $a<1$ and $d<1$.

Hence, we have shown that the system (\ref{eq-G31})-(\ref{eq-G32}) is
solved by (\ref{eq-sol24}), (\ref{eq-sol3}), and $s_1=\alpha\,s_3$, where
$\alpha$ is a root of $P_3$ satisfying (\ref{eq-P12}). In numerical
experiments we found that any root of $P_3$ satisfies (\ref{eq-P12}). Note
that since both tensors in (\ref{eq-SX2}) have multilinear rank $(2,2,2)$,
it follows that ${\bf S}$ is nonsingular. Hence, we have shown that for any
$\ul{\bf X}$ in orbit $G_3$ there exists a nonsingular ${\bf S}$ such that
$({\bf S},{\bf S},{\bf S})\cdot\ul{\bf Y}=\ul{\bf X}$, where $\ul{\bf Y}$
is the canonical form of orbit $G_3$.\\~\\

\noindent {\bf Acknowledgment.} The authors would like to thank Jos ten
Berge for commenting on an earlier version of this paper, and for drawing
their attention to $p\times p\times 2$ tensors.


\begin{thebibliography}{99}
\baselineskip14pt
\bibitem{BCS} {\sc P. B\"urgisser, M. Clausen and M.A. Shokrollahi} (1997)
{\em Algebraic Complexity Theory}, Springer, Berlin.

\bibitem{CC} {\sc J.D. Carroll and J.J. Chang} (1970) Analysis of
individual differences in multidimensional scaling via an $n$-way
generalization of Eckart-Young decomposition. {\em Psychometrika}, {\bf
35}, 283--319.

\bibitem{pierre} {\sc P. Comon} (1994) Independent component analysis, a
new concept? {\em Signal Processing}, {\bf 36}, 287--314.

\bibitem{ComMou} {\sc P. Comon and B. Mourrain} (1996) Decomposition of
quantics in sums of powers of linear forms. {\em Signal Processing}, {\bf
53}, 93--107.

\bibitem{Com98} {\sc P. Comon} (1998) Blind channel identification and
extraction of more sources than sensors. Keynote address at the {\em SPIE
Conference}, San Diego, July 19--24, pp. 2--13.

\bibitem{Comon} {\sc P. Comon} (2002) Tensor decompositions. In: {\em
Mathematics in Signal Processing}, Vol. V, J.G. McWhirter, I.K. Proudler
(Eds.), Clarendon Press, Oxford, USA.

\bibitem{CGLM} {\sc P. Comon, G. Golub, L.-H. Lim and B. Mourrain} (2008)
Symmetric tensors and symmetric tensor rank. {\em SIAM Journal on Matrix
Analysis and Applications}, {\bf 30}, 1254--1279.

\bibitem{AlmFM} {\sc A.L.F. de Almeida, G. Favier and J.C.M. Mota} (2007)
Parafac-based unified tensor modeling for wireless communication systems.
{\em Signal Processing}, {\bf 87}, 337--351.

\bibitem{LMVrank1} {\sc L. De Lathauwer, B. De Moor and J. Vandewalle}
(2000) On the best rank-1 and rank-$(R_1,R_2,\ldots,R_N)$ approxiamtion of
higher-order tensors. {\em SIAM Journal on Matrix Analysis and
Applications}, {\bf 21}, 1324--1342.

\bibitem{LMV2} {\sc L. De Lathauwer, B. De Moor and J. Vandewalle}
(2000). An introduction to independent component analysis. {\em Journal of
Chemometrics}, {\bf 14}, 123--149.

\bibitem{LC} {\sc L. De Lathauwer and J. Castaing} (2007) Tensor-based
techniques for the blind separation of DS-CDMA signals. {\em Signal
Processing}, {\bf 87}, 322--336.

\bibitem{DSL} {\sc V. De Silva and L.-H. Lim} (2008) Tensor rank and the
ill-posedness of the best low-rank approximation problem. {\em SIAM
Journal on Matrix Analysis and Applications}, {\bf 30}, 1084--1127.

%\bibitem{matcomp} {\sc G.H. Golub and C.F. Van Loan} (1996) {\em Matrix
%Computations, third edition}. Baltimore: John Hopkins University Press.

\bibitem{Har} {\sc R.A. Harshman} (1970) Foundations of the Parafac
procedure: models and conditions for an ``explanatory" multimodal factor
analysis. {\em UCLA Working Papers in Phonetics}, {\bf 16}, 1--84.

\bibitem{H1} {\sc F.L. Hitchcock} (1927) The expression of a tensor or a
polyadic as a sum of products, {\em Journal of Mathematics and Physics},
{\bf 6}, 164--189.

\bibitem{H2} {\sc F.L. Hitchcock} (1927) Multiple invariants and
generalized rank of a $p$-way matrix or tensor, {\em Journal of Mathematics
and Physics}, {\bf 7}, 39--70.

\bibitem{Hyv} {\sc A. Hyv\"arinen, J. Karhunen and E. Oja} (2001) {\em
Independent Component Analysis}, New York: Wiley.

\bibitem{JJ} {\sc J. Ja' Ja'} (1979) Optimal evaluation of pairs of
bilinear forms. {\em SIAM Journal on Computing}, {\bf 8}, 443--462.

\bibitem{KofReg} {\sc E. Kofidis and P.A. Regalia} (2002) On the best
rank-1 approximation of higher-order supersymmetric tensors. {\em SIAM
Journal on Matrix Analysis and Applications}, {\bf 23}, 863--884.

\bibitem{KolBad} {\sc T.G. Kolda and B.W. Bader} (2009) Tensor
decompositions and applications. {\em SIAM Review}, to appear.

\bibitem{KDS} {\sc W.P. Krijnen, T.K. Dijkstra and A. Stegeman} (2008) On
the non-existence of optimal solutions and the occurrence of ``degeneracy"
in the Candecomp/Parafac model. {\em Psychometrika}, {\bf 73}, 431--439.

\bibitem{Kro} {\sc P.M. Kroonenberg} (2008) {\em Applied Multiway Data
Analysis}, Wiley Series in Probability and Statistics.

\bibitem{K} {\sc J.B. Kruskal} (1977) Three-way arrays:
rank and uniqueness of trilinear decompositions, with applications to
arithmetic complexity and statistics. {\em Linear Algebra and its
Applications}, {\bf 18}, 95--138.

\bibitem{Krus89} {\sc J.B. Kruskal} (1989) Rank, decomposition, and
uniqueness for 3-way and $N$-way arrays, pp. 7--18 in: {\em Multiway Data
Analysis}, R. Coppi and S. Bolasco (Eds.), North-Holland.

\bibitem{KR84} {\sc J.P.S. Kung and G.-C. Rota} (1984) The invariant theory
of binary forms. {\em Bulletin of the American Mathematical Society}, {\bf
10}, 27--85.

%\bibitem{Lim-palo} {\sc L.-H. Lim} (2004) What's possible and what's
%impossible in tensor decompositions/approximation. Talk at the {\em Tensor
%Decomposition Workshop}, July 19--23, Palo Alto, USA.

%\bibitem{Lim2005} {\sc L.-H. Lim} (2005) Optimal solutions to non-negative
%Parafac/multilinear NMF always exist. Talk at the {\em Workshop on Tensor
%Decompositions and Applications}, August 29 - September 2, CIRM, Luminy,
%Marseille, France.

\bibitem{Ni} {\sc G. Ni and Y. Wang} (2007) On the best rank-1
approximation to higher-order symmetric tensors. {\em Mathematical and
Computer Modelling}, {\bf 46}, 1345--1352.

\bibitem{Qi} {\sc L. Qi} (2006) Rank and eigenvalues of a supersymmetric
tensor, the multivariate homogeneous polynomial and the algebraic
hypersurface it defines. {\em Journal of Symbolic Computation}, {\bf 41},
1309--1327.

\bibitem{SGB} {\sc N. Sidiropoulos, G. Giannakis and R. Bro} (2000) Blind
Parafac receivers for DS-CDMA systems. {\em IEEE Transactions on Signal
Processing}, {\bf 48}, 810--823.

\bibitem{SBG2} {\sc N. Sidiropoulos, R. Bro, G. Giannakis} (2000) Parallel
factor analysis in sensor array processing. {\em IEEE Transactions on Signal
Processing}, {\bf 48}, 2377--2388.

\bibitem{SBG} {\sc A. Smilde, R. Bro and P. Geladi} (2004) {\em Multi-way
Analysis: Applications in the Chemical Sciences}. Wiley.

\bibitem{Ste} {\sc A. Stegeman} (2006) Degeneracy in Candecomp/Parafac
explained for $p\times p\times 2$ arrays of rank $p+1$ or higher. {\em
Psychometrika}, {\bf 71}, 483--501.

\bibitem{Ste3} {\sc A. Stegeman} (2007) Degeneracy in Candecomp/Parafac
explained for several three-sliced arrays with a two-valued typical rank.
{\em Psychometrika}, {\bf 72}, 601--619.

\bibitem{Ste2} {\sc A. Stegeman} (2008) Low-rank approximation of generic
$p\times q\times 2$ arrays and diverging components in the
Candecomp/Parafac model. {\em SIAM Journal on Matrix Analysis and
Applications}, {\bf 30}, 988--1007.

\bibitem{SS} {\sc A. Stegeman and N.D. Sidiropoulos} (2007) On
Kruskal's uniqueness condition for the Candecomp/ Parafac decomposition.
{\em Linear Algebra and its Applications}, {\bf 420}, 540--552.

\bibitem{SDL} {\sc A. Stegeman and L. De Lathauwer} (2009) A method to avoid
diverging components in the Candecomp/Parafac model for generic $I\times
J\times 2$ arrays. {\em SIAM Journal on Matrix Analysis and
Applications}, {\bf 30}, 1614--1638.

\bibitem{Stras} {\sc V. Strassen} (1983) Rank and optimal computation of
generic tensors. {\em Linear Algebra and its Applications}, {\bf 52},
645--685.

\bibitem{BKL} {\sc J.M.F. Ten Berge, H.A.L. Kiers and J. De Leeuw} (1988)
Explicit Candecomp/Parafac solutions for a contrived $2\times 2\times 2$
array of rank three. {\em Psychometrika}, {\bf 53}, 579--584.

\bibitem{TBK} {\sc J.M.F. Ten Berge and H.A.L. Kiers} (1999) Simplicity of
core arrays in three-way principal component analysis and the typical rank
of $p\times q\times 2$ arrays. {\em Linear Algebra and its Applications},
{\bf 294}, 169--179.

\bibitem{TBSR} {\sc J.M.F. Ten Berge, N.D. Sidiropoulos and R. Rocci}
(2004) Typical rank and indscal dimensionality for symmetric three-way
arrays of order $I\times 2\times 2$ or $I\times 3\times 3$. {\em Linear
Algebra and its Applications}, {\bf 388}, 363--377.

%\bibitem{Tuk} {\sc L.R. Tucker} (1966) Some mathematical notes on
%three-mode factor analysis. {\em Psychometrika}, {\bf 31}, 279--311.

\bibitem{Zhang} {\sc T. Zhang and G. Golub} (2001) Rank-one approximation
to high order tensors. {\em SIAM Journal on Matrix Analysis and
Applications}, {\bf 23}, 534--550.

\end{thebibliography}
\end{document}